\documentclass[english, side, 10pt]{article}	
\usepackage[T1]{fontenc}			
\usepackage[utf8]{inputenc}			
\usepackage[verbose, a4paper]{geometry}		
\usepackage{color}					
\usepackage[english]{babel}			
\usepackage{array}
\usepackage{caption}				
\usepackage{subcaption}				
\usepackage{courier}
\usepackage{booktabs}
\usepackage{siunitx}
\usepackage{wrapfig} 				
\usepackage{amssymb, amsthm, amsmath, units}
\usepackage{tikz} 					
\usepackage{epstopdf}
\usepackage{todonotes}
\usepackage{float}
\usepackage[ruled,  vlined, linesnumbered]{algorithm2e} 

\usepackage{eurosym} 
\usepackage{enumitem}

\usepackage{newtxtext,newtxmath}	
\usepackage{dsfont} 				
\usepackage{microtype}   			

\usepackage[round]{natbib} 			
\usepackage[bookmarks=true, breaklinks= true, backref= page, 
pdftitle={Tree Generation}, pdfauthor={Georg Pflug, Alois Pichler, Kipngeno Kirui}, colorlinks=true, citecolor=blue, urlcolor=blue, pdfstartview= FitV]{hyperref}
\usepackage{doi} 					
\usepackage[capitalise]{cleveref}	

\theoremstyle{plain}
\newtheorem{theorem}       {Theorem}
\newtheorem*{theorem*}	   {Theorem}

\theoremstyle{definition}

\theoremstyle{remark}
\newtheorem*{remark*}{Remark}
\newtheorem{example}[theorem]{Example}
\newtheorem*{example*}{Example}

\DeclareMathOperator{\E}{\mathds E}				
\DeclareMathOperator{\var}{\mathsf{var}}

\newcommand{\pred}{\operatorname{pred}}

\DeclareMathOperator*{\argmin}{\operatornamewithlimits{arg\,min}}
\newcommand{\dd}{\mathsf {d\kern -0.07em l}} 
\newcommand{\nestr}{\dd} 
\newcommand{\Fc}{\mathcal{F}}

\begin{document}
\title{\textbf{New Algorithms And Fast Implementations To Approximate Stochastic Processes}}
\author{
	Kipngeno Benard Kirui\thanks{
		Technische Universität Chemnitz, Fakultät für Mathematik, Germany}
	\and
	Georg Ch.\ Pflug\thanks{
		University of Vienna and International Institute for Applied System Analysis (IIASA), Laxenburg, Austria}
	\and
	Alois Pichler\thanks{
		Contact: \protect\href{mailto:alois.pichler@math.tu-chemnitz.de}{alois.pichler@math.tu-chemnitz.de}.
		Deutsche Forschungsgemeinschaft (DFG, German Research Foundation)~-- Project-ID 416228727 -- SFB~1410. \protect\href{https://orcid.org/0000-0001-8876-2429}{https://orcid.org/0000-0001-8876-2429}.}~\,\footnotemark[1]
}
\maketitle

\begin{abstract}
We present new algorithms and fast implementations to find efficient approximations for modelling stochastic processes. 
For many numerical computations it is essential to develop finite approximations for stochastic processes. 
While the goal is always to find a finite model, which represents a given knowledge about the real data process as accurate as possible, the ways of estimating the discrete approximating model may be quite different: (i)~if the stochastic model is known as a solution of a stochastic differential equation, e.g., one may generate the scenario tree directly from the specified model; (ii)~if a simulation algorithm is available, which allows simulating trajectories from all conditional distributions, a scenario tree can be generated by stochastic approximation; (iii)~if only some observed trajectories of the scenario process are available, the construction of the approximating process can be based on non-parametric conditional density estimates.

We also elaborate on the important concept of distances, which allows us to assess the quality of the approximation. We study these methods and apply them to an electricity price data.

Our fast implementation
\begin{center}
	\href{https://kirui93.github.io/ScenTrees.jl/stable/}{\texttt{ScenTrees.jl}}
\end{center} including an exhaustive documentation is available for free at GitHub.\footnote{\label{fn:ScenTrees}\texttt{ScenTrees.jl}: \href{https://kirui93.github.io/ScenTrees.jl/stable/}{https://kirui93.github.io/ScenTrees.jl/stable/}, cf.\ \citet{KiruiPichler}}
\medskip


	
\textbf{Classification:} 90C15, 60B05, 62P05
\end{abstract}

\section{Introduction}
	A decision maker can formulate an action plan to achieve a specific goal. 
	The action plan is executed in sequential stages.
	Planning under uncertainty considers already at the beginning random events which might happen while the plan is executed, taking the possibility of corrective (recourse) decisions into account. This ensures that the desired goal can be met, even in an uncertain environment where random events impact and influence the situation at each stage.

	An important historical contribution in the development of decision making under uncertainty is \citet{Markowitz1952}, who considers an optimal asset allocation problem and developed what is nowadays called \emph{modern portfolio theory}. This example is often used to illustrate and outline the difficulty: an investment decision has to be made today, but the random return is observed after the investment horizon, in the future. In its simplest version, the problem has only one decision stage, but if a recourse decision is allowed after the random outcome is revealed, the problem is said to be \emph{two-stage}. However, if financial securities such as options with different maturities are included in the portfolio then more decision stages must be considered: intermediate decisions have to be made whenever options expire before reaching the investment horizon; the problem is then said to be \emph{multistage}.


	Further problems dealing with decision making under uncertainty appear in logistics and inventory control. Goods must be shipped between various warehouses in view of random demands and unmatched demand triggers recourse actions by rapid new orders (cf.\ \citet{Dantzig1955} for a very early problem description).
	This problem can be formulated as a decision problem in multiple stages. Optimal inventory and transportation policies constitute hedging strategies against unknown and varying demands.

	As another problem with economic relevance we mention an optimal schedule to manage reservoirs of a hydropower plant. Its inflows are rain and perhaps melting snow, which are both random events. The price of electricity produced is highly volatile in addition, so that the optimal operation mode is a sequence of decisions, made at consecutive times and influenced by various random events (cf.\ \citet{Seguin2015} for a more concrete discussion).
	\medskip

	All problems mentioned above have uncertain scenarios and decisions at intermediate stages in common. The scenarios appear as the trajectories of a stochastic process in discrete time, $\xi=(\xi_1, \dots, \xi_T)$, with values $\xi_t\in\mathbb{R}^{m_t}$, say. For algorithmic treatment, such general processes must be approximated by simpler ones, in the same manner as functions are represented as vectors on finite grids on digital computers. In particular, a discrete time stochastic process is approximated by a \emph{scenario tree} or, in case of a Markovian process, by a \emph{scenario lattice}. Notice, however, that quite often the decisions at earlier stages influence the feasibility of profitability of later stages and in such a case the full history of the decision process has to be recorded. The history process is always a tree.

	\medskip
	An elementary component for generating scenario trees is the approximation of a random outcome, i.e., the approximation of random variables. An example is the random realization~$\xi_1$ at the first stage, which has to be approximated sufficiently well for an algorithmic treatment.
	\citet{GrafLuschgy} provide a comprehensive study and theoretical background on this topic (called also {\it uncertainty vector quantification}), \citet{PagesFunctionalQuantization} and \citet{PagesBallyPrintems} present methods and algorithms, many of them are based on optimal clustering, cf.\ \citet{Hartigan1975}. It is an important observation that the approximation of a random outcome as $\xi_1$ suffers from the curse of dimensionality, a term coined by R.\ Bellman (cf.\ \citet{Bellman1961,Dudley1969}).

	The generation of scenario trees is even more involved. Indeed, not only one random realization as  $\xi_1$ has to be approximated, but the entire underlying stochastic process $(X_1,\dots,X_T)$, often with complicated interdependencies. Consequently, the approximation of the stochastic process is significantly more complex and computationally more expensive.

	\medskip \citet{Hoyland2001} describe an early method to generate scenario trees, they employ selected \emph{statistical parameters} such as moments and correlations to characterize transitions and interdependences within a tree. Their method 
	requires solving a system of nonlinear equations to find representative discretizations of the scenario tree. However, the corresponding transition probabilities ignore all further indispensable characteristics of the governing distribution and for that the method is not appropriate in general. \citet{Klaassen2002} discusses the moment matching method in relation to arbitrage opportunities. In a broader context, the method is also outlined in \citet{CONSIGLI2000} including an algorithm employing sequential updates and a concrete managerial example. \citet{Kuhn2015} address an asset allocation problem in an ambiguous context.

	Some decision problems exhibit specific problem characteristics. These properties are occasionally addressed directly by adjusting and designing the scenarios adequately, as \citet{HenrionRomisch2017} propose. The recent work \citet{Wallace2015}, e.g., particularly puts the tails of the distributions in focus, while \citet{KautWallace2011} involve copulas to properly capture the multivariate shape of the distributions.

	\cite{Beyer2014} develop an evolution strategy to generate scenario trees. They employ swarm intelligence techniques to tackle related multistage optimization problems for optimal portfolio decisions.
	\citet{Keutchayan2017} present a method to generate large scale scenario trees. They observe that the scenario tree has to be reduced for computational efficiency. The paper presents a heuristic which reduces the loss of information which comes along with reducing the scenario tree. \citet{HeitschRomisch2009} propose advanced heuristics as well to reduce large scale scenario to a computationally tractable size.

	The quality of a scenario trees can be assessed if there is a distance available, which measures the distance to the genuine stochastic process to be approximated. \cite{Pflug2009} presents the first concept to measure the distance of different stochastic processes. \citet{KovacevicPichler} employ the method directly to generate scenario trees.

	\medskip
	This paper endorses nonparametric techniques for generating scenario trees. In contrast to most methods addressed above, they involve the full law (distribution) of the stochastic process rather than individual, selected statistical parameters. In this way the marginals, moments, covariance and other characteristics are correctly mirrored and asymptotic consistency can be shown. The methods do not require specific parametric assumptions on the underlying process and covers multivariate observations as well. To construct an approximating tree it is further enough to have samples from the process available, these samples are just trajectories or sample paths. We build the trees by employing stochastic approximation (see \citet{Kushner2008} for a survey). It is not necessary to specify and solve a system of corresponding equations.

\paragraph{Implementation.}
Our Julia implementation is available freely on GitHub in the package \texttt{ScenTrees.jl},\footnote{\texttt{ScenTrees.jl}: \href{https://github.com/kirui93/ScenTrees.jl}{https://github.com/kirui93/ScenTrees.jl}} which is released under the open-source MIT license. Description of various functions of the package and also different examples to illustrate the different methods can be found in the package's documentation.\footnote{\texttt{Documentation}: \href{https://kirui93.github.io/ScenTrees.jl/stable/}{https://kirui93.github.io/ScenTrees.jl/stable/}, cf.\ \citet{KiruiPichler}}

\paragraph{Structure of the paper.}
In~\cref{sec:Notation}, we introduce some prerequisites that are necessary for the discussion in this paper. This setting includes a discussion on the important concept of the quality of approximation in \cref{sec:Nested}. 
\cref{sec:Trajectories} outlines the methods of scenario tree generation by nested clustering.
Successive improvements are obtained by adding further information or trajectories to a tree and \cref{sec:StochApproximation} presents the stochastic approximation procedure.
\cref{sec:Trees} addresses the situation of limited data available, as is the case in many situations of practical interest.
\cref{sec:Lattice} is dedicated specifically to stochastic processes with special properties as Markovian processes or processes with fixed states.
A brief introduction on trees without predefined structure is addressed in \cref{sec:Treeswithout}. Lastly, \cref{sec:application} presents a special application of scenario lattice generation methods on a limited electricity load data. We draw our conclusions on \cref{sec:Summary}.

\section{Mathematical setting}\label{sec:Notation}
In what follows we consider a general stochastic process $X$ in discrete time over $T$ stages, $X = (X_1,\ldots,X_T)$, where $X_1$ is a nonrandom starting value.
A scenario tree is a discrete time and discrete state process approximating the process~$X$. Throughout this paper, we employ the following terminology and conventions for scenario trees, scenario lattices and observed data. 

\subsection{Scenario tree}
A finite scenario tree approximating~$X$ is denoted by $\tilde X=\left(\tilde X_1, \tilde X_2,\dots,\tilde X_T\right)$ and modelled as a scenario tree process.

\begin{figure}[H]	
	\centering
	\begin{tikzpicture}[scale=0.7]
	\node (1)  at (0,0) [circle,shade,draw] {1};
	\node (2)  at (2,2)[circle,shade,draw] {2};
	\node (3)  at (2,0) [circle,shade,draw] {3};
	\node (4)  at (2,-1.5)[circle,shade,draw] {4};
	\node (5)  at (4,3.1) [circle,shade,draw] {5};
	\node (6)  at (4,2.1) [circle,shade,draw] {6};
	\node (7)  at (4,1.1) [circle,shade,draw] {7};
	\node (8)  at (4,0) [circle,shade,draw] {8};
	\node (9)  at (4,-1)[circle,shade,draw] {9};
	\node (10) at (4,-2)[circle,shade,draw] {10};
	
	\node (13) at (8,-2.5)[circle,shade,draw] {$\mathbf n$};
	\node (12) at (8,-1.5)[circle,shade,draw] {...};
	\node (14) at (8,3.5)[circle,shade,draw] {...};
	\node (15) at (8,2.5)[circle,shade,draw] {...};
	\node (16) at (8,0.5)[circle,shade,draw] {...};
	\node (17) at (8,-.5)[circle,shade,draw] {...};
	
	\draw[->, very thick] (1) to (2);
	\draw[->, very thick] (2) to (5);
	\draw[->, very thick] (2) to (6);
	\draw[->, very thick] (2) to (7);
	\draw[->, very thick] (1) to (3);
	\draw[->, very thick] (3) to (8);
	\draw[->, very thick] (1) to (4);
	\draw[->, very thick] (4) to (9);
	\draw[->, very thick] (4) to (10);
	
	\draw[->, dotted, very thick] (5,-2) to (5.7,-2);
	\draw[->, dotted, very thick] (6,-2) to (13);
	\draw[->, dotted, very thick] (6,-2) to (12);
	
	\draw[->, dotted, very thick] (5,0) to (5.7,0);
	\draw[->, dotted, very thick] (6,0) to (16);
	\draw[->, dotted, very thick] (6,0) to (17);
	
	\draw[->, dotted, very thick] (5,3) to (5.7,3);
	\draw[->, dotted, very thick] (6,3) to (14);
	\draw[->, dotted, very thick] (6,3) to (15);

	\node (t1) at (0,-3.5) {$t=1$};
	\node (t2) at (2,-3.5) {$t=2$};
	\node (t3) at (4,-3.5) {$t=3$};
	\node (t4) at (8,-3.5) {$t=T$};
	\end{tikzpicture}
	\caption{Scenario tree process in $T$ stages\label{fig:Tree}}
\end{figure}
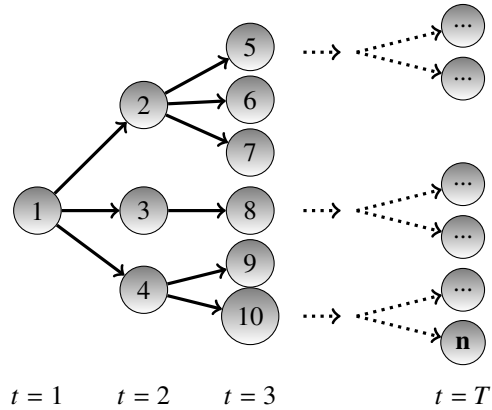

\paragraph{Characterization of a scenario tree process.}
A directed, rooted and layered graph describes the topology of the scenario tree.
We enumerate all $n$ nodes of a scenario tree by $1,\dots,n$, the root node of the scenario tree is $1$ (the nodes $\mathcal N:=\{1,\dots,n\}$ are the vertices of the graph, cf.\ \cref{fig:Tree}).
A scenario tree process is then fully characterized by the following three lists:
\begin{enumerate}[label=(\roman*)]
	\item The list of predecessors,
	\[\pred(1),\dots,\pred(n),\]
	determines the \emph{topology} of the tree. $\pred(i)$ is the predecessor of node number~$i$. For the root node we set $\pred(1):=0$. (For the exemplary tree in \cref{fig:Tree}, for instance, $\pred(2)=1$ and $\pred(8)=3$.) The pairs $(\pred(i),i)$, $i=2,\dots,n$, comprise all arcs of the graph.

	\item The list of \emph{probabilities} is
	\[\tilde p_1,\dots, \tilde p_n,\]
	where $\tilde p_i$ denotes the conditional probability to reach the node $i$ from its predecessor along the edge of the tree. We set $\tilde p_1=1$ for the root node.
	
	\item The list of \emph{states} (i.e., the values of the tree process at the nodes) is \[\tilde x_1,\dots,\tilde x_n.\]
	The values $\tilde x_i$ constitute the outcomes of the process $\tilde X$ at stage $t$, they may be vectors of any dimension. We denote the dimension of the process at stage $t$ by $m_t$, i.e., $\tilde x_i\in\mathbb R^{m_t}=:\Xi_t$, whenever the node $i$ is at stage~$t$.\footnote{To allow a compact presentation below we shall assume, without loss of generality, that $\Xi_t=:\Xi$ for all $t\le T$.}
\end{enumerate}

A sample scenario is a path $(j_1,\dots,j_T)$ with  
$j_{t-1}=\pred(j_t)$. 
The sample scenario connects the root $1$ and a leaf $j_T$. A scenario tree has as many scenarios as there are leaves. The probability of a scenario is $P\big((j_1,\dots,j_T)\big)=\tilde{p}_{j_1}\cdot\ldots\cdot \tilde{p}_{j_T}$.

The predecessor of the node $j\in\mathcal N$ at stage $t$ is $j_t$, that is, there is a sequence $(j_1,\dots,j_t,\dots, j_\tau)$ for some $\tau\ge t$ with $j_\tau=j$ and $\pred(j_{t^\prime})=j_{t^\prime-1}$ for all $t^\prime\le\tau$. (For the exemplary tree in \cref{fig:Tree}, for example, the predecessor of the node $n$ at stage $2$ is $4$, i.e., $n_2=4$.)


Filtrations model the flow of information associated with the tree process. The filtration is the sequence of $\sigma$-algebras $\mathcal F_t$ modeling the information available at stage $t$. The $\sigma$-algebra $\mathcal F_t$ is generated by
\[\mathcal F_t=\sigma\big( \left\{s=(j_1,\dots,j_T)\colon s\text{ is a scenario with } j_t=j\right\}, j\in\mathcal N\big).\]
We mention that the $\sigma$-algebras $\mathcal F_t$, $t=1,\dots,T$, uniquely characterize the topology, i.e., the branching structure of the process.  Notice that $\tilde{X}$ is a tree process, if the $\sigma$-algebra generated by $\tilde{X}_t$ is the same as the one generated by $(\tilde{X}_1, \dots, \tilde{X}_t)$ for all $t$.

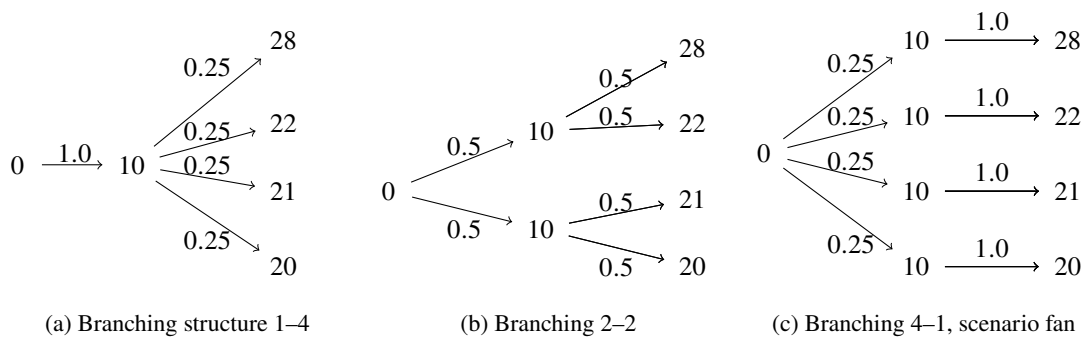
\begin{figure}[H]
	\centering
	\unitlength0.5cm
	\hfill
	\begin{subfigure}[t]{0.33\textwidth}
		\begin{tikzpicture}[scale=0.5]
		\node (1)  at (0,0) [circle] {0};
		\node (2)  at (3.0,0) [circle] {10};
		\node (3)  at (7,1.1)[circle] {22};
		\node (4)  at (7,-0.7) [circle] {21};
		\node (5)  at (7,3.3) [circle] {28};
		\node (6)  at (7,-2.7) [circle] {20};
		
		\foreach \x in {2} 	{
			\draw[->,thin] (1) to (\x);
			\foreach \y in {3,4,5,6} \draw[->, thin] (\x) to (\y);}
		\draw (1.5,0.3) node{$1.0$};
		\draw (5.0,2.6) node{$0.25$};
		\draw (5.0,0.9) node{$0.25$};
		\draw (5.0,0.0) node{$0.25$};
		\draw (5.0,-2.0) node{$0.25$};
		\end{tikzpicture}
		\subcaption{Branching structure 1--4}\label{fig:t1}
	\end{subfigure}\hfill
	\begin{subfigure}[t]{0.33\textwidth}
		\begin{tikzpicture}[scale=0.5]
		\node (1)  at (0,0) [circle] {0};
		\node (2)  at (4,1.6) [circle] {10};
		\node (3)  at (4,-1.0)[circle] {10};
		\node (4)  at (8,3.8) [circle] {28};
		\node (5)  at (8,-0.2) [circle] {21};
		\node (6)  at (8,1.8) [circle] {22};
		\node (7)  at (8,-2.0) [circle] {20};
		
		\foreach \x in {2,3}{\draw[->,thin] (1) to (\x);
			\foreach \y in {5,7} \draw[->, thin] (3) to (\y);
			\foreach \y in {4,6} \draw[->, thin] (2) to (\y);}
		\draw (2,1.2) node{$0.5$};
		\draw (2,-1.0) node{$0.5$};
		\draw (6,3.0) node{$0.5$};
		\draw (6,2.0) node{$0.5$};
		\draw (6,-0.3) node{$0.5$};
		\draw (6,-2.0) node{$0.5$};
		\end{tikzpicture}
		\subcaption{Branching 2--2}\label{fig:t2}
	\end{subfigure}
	\begin{subfigure}[t]{0.33\textwidth}
		\begin{tikzpicture}[scale=0.5]
		\node (1)  at (0,0) [circle] {0}; 
		\node (2)  at (4,1.0) [circle] {10};
		\node (3)  at (4,-1.0)[circle] {10};
		\node (8)  at (4,3.0)[circle] {10};
		\node (9)  at (4,-3.0)[circle] {10};
		\node (4)  at (8,1.0) [circle] {22};
		\node (5)  at (8,-1.0) [circle] {21};
		\node (6)  at (8,3.0) [circle] {28};
		\node (7)  at (8,-3.0) [circle] {20};
		
		\foreach \x in {2,3,8,9}{\draw[->,thin] (1) to (\x);
			\foreach \y in {5} \draw[->, thin] (3) to (\y);
			\foreach \y in {7} \draw[->, thin] (9) to (\y);
			\foreach \y in {4} \draw[->, thin] (2) to (\y);
			\foreach \y in {6} \draw[->, thin] (8) to (\y);}
		\draw (2.3,1.0) node{$0.25$};
		\draw (2.3,2.4) node{$0.25$};
		\draw (2.3,-0.2) node{$0.25$};
		\draw (2.3,-2.4) node{$0.25$};
		\draw (6,3.5) node{$1.0$};
		\draw (6,1.5) node{$1.0$};
		\draw (6,-0.5) node{$1.0$};
		\draw (6,-2.5) node{$1.0$};
		\end{tikzpicture}
		\subcaption{Branching 4--1, scenario fan}\label{fig:t3}
	\end{subfigure}\hfill
	\caption{Three scenario trees with identical paths but still different topologies}
	\label{fig:1}
\end{figure}

\medskip
\cref{fig:1} presents three scenario tree processes with different topology and filtrations. Nonetheless, the trajectories of their states are identical, their scenarios even have the same probabilities.
\subsection{Scenario lattices}
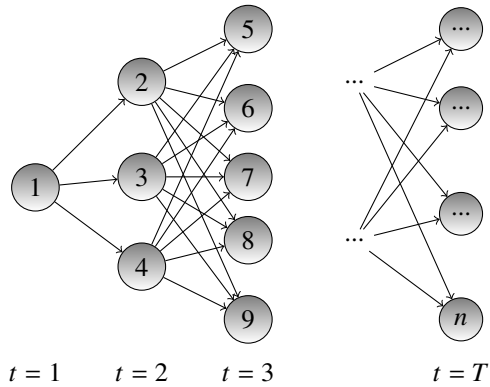
\begin{figure}[H]	
	\centering
	\begin{tikzpicture}[scale= 0.7]
	\node (1)  at (0,0) [circle,shade,draw] {1};
	\node (2)  at (2,2)[circle,shade,draw] {2};
	\node (3)  at (2,0.2) [circle,shade,draw] {3};
	\node (4)  at (2,-1.5)[circle,shade,draw] {4};
	\node (5)  at (4,3) [circle,shade,draw] {5};
	\node (6)  at (4,1.5) [circle,shade,draw] {6};
	\node (7)  at (4,0.2) [circle,shade,draw] {7};
	\node (8)  at (4,-1.0) [circle,shade,draw] {8};
	\node (9)  at (4,-2.5)[circle,shade,draw] {9};
	
	\node (13) at (8,-2.5)[circle,shade,draw] {$n$};
	\node (12) at (8,-0.5)[circle,shade,draw] {...};
	\node (14) at (8,3)[circle,shade,draw] {...};
	\node (15) at (8,1.5)[circle,shade,draw] {...};

	\node (A) at (6,2) {...};
	\node (B) at (6,-1) {...};
	
	\foreach \x in {2, 3, 4} 	{
		\draw[->,thin] (1) to (\x);
		\foreach \y in {5,6,7,8,9} \draw[->, thin] (\x) to (\y);}
	\foreach \x in {12,...,15} 	{
		\draw[->,thin] (A) to (\x); \draw[->, thin] (B) to (\x);}
	
	\node (t1) at (0,-3.5) {$t=1$};
	\node (t2) at (2,-3.5) {$t=2$};
	\node (t3) at (4,-3.5) {$t=3$};
	\node (t4) at (8,-3.5) {$t=T$};
	\end{tikzpicture}
	\caption{A model of a lattice at $T$ stages and $n$ nodes\label{fig:Lattice}}
\end{figure}

Scenario lattices are a natural discretization of Markovian processes. While nodes of a scenario tree have one predecessor, nodes of a scenario lattice can have multiple predecessors as shown in \cref{fig:Lattice}. Many stochastic processes which occur in real world applications represent a Markovian process, i.e., the further evolution of the process depends on the state at the actual stage, but not on the history (the process is colloquially said to be \emph{memoryless}). This additional knowledge can and should be explored in implementations, the consequences are manifold.

In the Markovian case it is often \emph{not} necessary to build a complete tree for optimization purposes. Instead, a \emph{lattice} (i.e., a recombining tree) may be sufficient to model the whole process, cf.\ \cref{fig:Lattice}. The lattice is the structure of a finitely valued Markovian process, which is a rooted, layered, directed graph of height $T$ (all paths form the root to the leaves have same length $T$), such that edges can only exist between nodes of subsequent layers.
\begin{remark*}
	The lattice model is not sufficient, e.g., if the decisions are path-dependent, which may happen even if the scenario model is Markovian. To give an example, suppose that the decisions contain to exercise  American or Asian options or to deal with knock-out certificates. In these cases, the whole history of the price process and/\,or the decision process is relevant at each state of the process and a scenario tree has to be considered anyway.
\end{remark*}

\subsection{The decision problem}
Once a scenario tree or lattice is available which models the underlying process properly, then the tree can be used for decision making to meet a predefined managerial goal. The goal to be met, e.g., is to maximize the discounted, expected profit from today's perspective, at the beginning of the planning horizon.

To achieve this goal, decisions have to be made at each node within the tree. This tree has the same structure as the scenario tree, it is called the \emph{decision tree}. As \citet{RockafellarWets1991} elaborate, the decisions made on the decision tree reflect a hedging strategy against the unknowns. The purpose of the hedging strategy is to meet the initial, predefined goal. 

\subsection{Quality of approximation}\label{sec:Nested}
The genuine process $X$ is often not eligible for decision making, this is the case if the process is available only from samples or because its structure is too complicated. An adequate approximation~$\tilde X$ of~$X$ should lead to reliable results, if decisions are based on $\tilde X$ instead of the initial process $X$.

The following subsections address the quality of approximation in two ways. First, we observe that the approximations $\tilde X$ obtained in the initial sections result by mapping the genuine process~$X$ to its tree approximation $\tilde X$, while the following section addresses the nested distance. Both perceptions make a distance available, which allows assessing the approximation quality in a quantitative way. More general results ensure that the policies, which are derived from the approximating process $\tilde X$, apply for the initial process $X$ as well and the model error can be stated explicitly in terms of the distances.

\subsubsection{Transportation maps}
Important procedures outlined below have in common that they fit individual trajectories into an existing tree structure.
This operation is a mapping, where trajectories are mapped on a tree by respecting the evolution of information over time (processes respecting the time evolution are said to be \emph{adapted}, or \emph{nonanticipative}).

To this end consider a map $\mathbf T$, mapping the state  $\xi=(\xi_1,\dots,\xi_T)$ of the stochastic process $X=(X_1,\dots X_T)$ to the new scenario $\left(\tilde x_1,\dots\tilde x_T\right):=\mathbf T(\xi)$. 
By the Doob--Dynkin lemma (cf.\ \citet{Kallenberg2002Foundations}), the \emph{transport map} $\mathbf T$ is adapted (i.e., it preserves the information) if its components are given by
\begin{equation*}\label{eq:T}
	\mathbf T\big(\xi_1,\dots,\xi_T\big)=
	\begin{pmatrix}\begin{array}{l}
		\tilde x_1=\mathbf T_1(\xi_1)\\
		\tilde x_2=\mathbf T_2(\xi_1,\xi_2)\\
		\quad\vdots\\
		\tilde x_T=\mathbf T_T(\xi_1,\dots,\xi_T)
	\end{array}\end{pmatrix}.
\end{equation*}

The trajectory~$\xi$ and the path $\tilde {\mathbf x}=\mathbf T(\xi)$ within the tree have the distance
\begin{equation}\label{eq:dist}
	d(\xi, \tilde{\mathbf x}):=\sum_{t=1}^T d_t\big(\xi_t,\tilde{\mathbf x}_t\big),
\end{equation}
where $d_t$ is the distance function on $\Xi_t$ for the  stage $t$. This allows defining and assessing the quality of the approximation by
\begin{equation}\label{eq:distT}
	\E_P\,d\big(\xi, \mathbf T(\xi)\big)= \int_{\Xi_1\times\dots\times\Xi_T} d\big(\xi,\,\mathbf T(\xi)\big)\,P(\mathrm d \xi),
\end{equation}
where the expectation (integration) is with respect to the law of the stochastic process~$X$.
The quantity~\eqref{eq:distT} is the \emph{average aberration} of the  initial stochastic process $X$ in comparison with the tree approximation specified by the mapping $\mathbf T$. Eq.~\eqref{eq:distT} is a reliable quantity describing the quality of the approximation provided by $\mathbf T$. The objective of usual stochastic optimization problems can be shown to depend continuously on the quantity~\eqref{eq:distT}. 
\medskip

To raise a conceptual issue associated with transport maps, consider an \emph{invertible} transport map $\mathbf T$ so that
\[	\E_P d\big(\xi,\mathbf T(\xi)\big)= \E_{P^{\mathbf T}} d\Big(\mathbf T^{-1}(\eta), \eta \Big)=\E_{P^{\mathbf T}} d\Big(\eta, \mathbf T^{-1}(\eta) \Big)\]
by the change of variable formula, where $P^{\mathbf T}:=P\circ \mathbf T^{-1}$ is the image measure.
A comparison with~\eqref{eq:distT} demonstrates that the inverse $\mathbf T^{-1}$ has to be adapted too to allow switching from the initial space to the image space.
This is, however, not the case in real-world application where $\mathbf T$ maps a stochastic process with continuous states to a tree process, which has only finitely many states.


The nested distance, presented in what follows, resolves this problem conceptually.

\subsubsection{Relation to the nested distance}
The \emph{nested distance} (or \emph{process distance}) is employed to measure the distance of two distinct stochastic processes directly. It is a symmetric version of~\eqref{eq:distT} and based on transportation distances derived from transportation theory, cf.\ \citet{Villani2003}.

The expression~\eqref{eq:distT} gives rise to rewrite the expectation as
\begin{equation}\label{eq:17}
	\E_P d\big(\xi,\mathbf T(\xi)\big)=\E_\pi d(\xi,\eta)=\iint  d(\xi,\eta)\,\pi(\mathrm d \xi,\mathrm d\eta),
\end{equation} employing the transport measure $\pi(A\times B):=P\big(A\cap\mathbf T^{-1}(B)\big)$.
As a generalization we consider general transportation measures $\pi$ in~\eqref{eq:17}, which are not necessarily of the form $P(A\cap\mathbf T^{-1}(B))$.
A general transport measure~$\pi$ is often termed \emph{transport plan} to distinguish it from the transportation map~$\mathbf T$.


This distance of probability measures obtained in this way is nowadays known as \emph{Kantorovich} or \emph{Wasserstein distance}, as \emph{earth mover}, or more generally as \emph{transportation distance}. We remark here that the notion of transportation distance extends obviously to non-discrete probability measures $P$ and $Q$ on a metric space $(\Xi,d)$.

The nested distance (see Appendix~\ref{sec:NestedA}) generalizes the transportation distance by additionally taking the increasing information at successive stages into account. In this way the nested distance generalizes transportation distances to a distance of stochastic processes, or to a distance of trees.

\section{Generation of scenario trees}\label{sec:Trajectories}
This section introduces various techniques of generating scenario trees. We also give examples of two stochastic processes that can be approximated via a scenario tree. 

\subsection{Data and observations}
The basic data for tree generation is typically given by a set of independent trajectories or sample paths, these trajectories constitute ideally observed data. Time series of stocks or energy prices, demands or reservoir inflow data, to provide examples, can be collected throughout comparable time intervals (weeks, say) and these data represent a reliable collection of trajectories. Alternatively, an econometric model could also be estimated and representative sample paths can be generated by simulation to serve as basic data.

Observations of the process $X$ are trajectories, denoted by
\begin{equation}\label{eq:2}\xi_j, \qquad j=1,\dots,N,\end{equation} with $N\in\mathbb N$ (or $N=\infty$, if the samples are generated by an econometric model). Each observation is a complete outcome of the process over time, i.e., $\xi_j=(\xi_{j,1},\dots,\xi_{j,T})$.

\subsection{Description of the methods}
The sample paths are employed to construct the scenario tree. In order to learn as much as possible from the data, the number of trajectories observed ($N$, cf.~\eqref{eq:2}) should be huge and representative in terms of spanning the underlying stochastic phenomena.
A new independent trajectory may be used at every new step of our algorithm, if enough data are available or may be generated by simulation.
The corresponding algorithms are presented below first. In some real life situations, however, the sample size may not be large and that is why we discuss later also methods involving kernel estimates, which work with limited data and small sample sizes.

The algorithms we propose start with a tree, which is not more than a qualified guess (based on an expert opinion, e.g.). Drawing one sample after another we improve the initial scenario tree successively using a technique which derives from stochastic approximation.
Each iteration modifies the values of the tree without changing its structure. In this way the approximating quality of the tree improves gradually. The tree fluctuates at the beginning, but with more and more iterations the tree converges in probability. The resulting tree finally is optimal with the properties desired and is ready to be used for the decision making procedure.

We also discuss a construction algorithm which starts with a small tree and gradually adds new branches when necessary. Using this technique one does not even need to have an initial guess of the tree structure.

In this paper we consider several practical situations which require algorithms for scenario tree generation:
\begin{itemize}[noitemsep]
	\item scenario tree generation from given stochastic models;
	\item scenario tree generation from stochastic simulations;
	\item scenario tree generation from observed  trajectories.
\end{itemize}

\subsection{Nested clustering}
The \emph{nested clustering} algorithm, which we outline in what follows, generates a tree from a \emph{finite} sample of $N$ trajectories, $\xi_1,\dots,\xi_N$. A finite sample of trajectories is occasionally called a scenario \emph{fan}, reminding to the fact that the slats of a fan do not bundle in their nodes (see, e.g., \cref{fig:t3}). The nested clustering algorithm determines the nodes of the approximation tree process~$\tilde X$ sequentially by iterating over stages, starting at stage $t=1$ up to the terminal stage $t=T$.

The initial state $\tilde X_1$ at stage $t=1$ of the process is assumed to be known (i.e., deterministic), $\tilde X_1= \tilde x_1$. The observations at the following stages $t=2,\dots, T$ are random.

The subsequent stage $t=2$ can be handled well. Note that the random component~$X_2$ at stage $t=2$ follows some distribution, which is known by assumption or from which data are available. By employing clustering techniques it is possible to approximate the distribution~$X_2$ by representative points, or cluster means. The $k$-means clustering algorithm, e.g., finds cluster means $\tilde x_i$ ($i\in\mathcal S_1$) on the stage $t=2$ which minimize the average minimal distance,\footnote{One can avoid double counting in~\eqref{eq:kMeans} of $\xi_{j,2}$, equidistant to various means $\tilde x_i$, by counting the distance only once, for example for the smallest node index $i$. We shall refer to this rule as the \emph{tie-break rule}.\label{fn:tie-break}}
\begin{equation}\label{eq:kMeans}
	\frac 1 N\sum_{j=1}^N \min_{i\in \mathcal S_1} d_2\left(\xi_{j,2},\, \tilde x_i\right),
\end{equation} where $d_2$ is the distance of outcomes at stage $t=2$ and \[\mathcal S_i:= \left\{\ell\in\mathcal N\colon\pred(\ell)=i\right\}\]
collects all successor-nodes of the node~$i$ (in \cref{fig:Tree}, $\mathcal S_1=\{2,3,4\}$ and $\mathcal S_2=\{5,6,7\}$, e.g.).
In this way,~$X_2$ is approximated by a discrete distribution $\tilde X_2$ located at $\tilde x_i$ ($i\in \mathcal S_1$), where each $\tilde x_i$ is reached with probability $P(\tilde X_2=\tilde x_i)=\tilde p_i$. 
The probabilities $\tilde p_i$ are found by allocating each sample $\xi_{j,2}$, $j=1,\dots N$, to the closest of all $\tilde x_i$ at the first stage, $i\in \mathcal S_1$. The probability $\tilde p_i$ thus represents the relative count of how often $\tilde x_i$ is the mean closest to $\xi_{j,2}$ (see the tie-breaking rule in \cref{fn:tie-break}). 

Note that the origin of the approximating tree is already found: the first two components of the approximating tree are $(\tilde X_1,\tilde X_2)$ with evaluations $(\tilde x_1, \tilde x_i)$ ($i\in \mathcal S_1$) at the first two stages.

\begin{figure}[H]	
	\centering
	\begin{tikzpicture}[scale= 1]
	\node (1)  at (0,0.2) {$\tilde{\boldsymbol x}_1$};
	\node (2)  at (2,2)  {$\tilde{\boldsymbol x}_2$};
	\node (3)  at (2,-0.0) {$\tilde{\boldsymbol x}_3$};
	\node (4)  at (2,-1.5){$\tilde{\boldsymbol x}_4$};
	\node (5)  at (4,3) {$\tilde{\boldsymbol x}_5$};
	\node (6)  at (4,1.9) {$\tilde{\boldsymbol x}_6$};
	\node (7)  at (4,.4) {$\tilde{\boldsymbol x}_7$};
	\node (8)  at (4,-0.5) {};
	\node (9)  at (4,0) {};
	\node (10) at (4,-1.5) {};
	
	\draw[->, ultra thick] (1) to (2);
	\draw[->, ultra thick] (2) to (5);
	\draw[->, ultra thick] (2) to (6);
	\draw[->, ultra thick] (2) to (7);
	\draw[->, dashed, very thin] (1) -- (2,-.3) --(4,1)--(4.3,0.9);
	\draw[->, dashed, very thin] (1) -- (2, .3)-- (4,0.1)--(4.3,0.3);
	\draw[->, dashed, very thin] (1) to (4) -- (4,-1.6)--(4.3,-1.5);
	\draw[->, dashed, very thin] (1) -- (2,-1.)-- (4,-.7)--(4.3,-1.0);
	
	\draw (2,0.3) --(2,1.7);
	\draw (2,2.3) --(2,3.5);
	\draw (2,-.2) --(2,-1.1);
	\draw (2,-1.8) --(2,-2.1);
	\draw[ultra thick] (1.9,-.7) --(2.1,-.7);
	\draw[ultra thick] (1.9,.8) --(2.1,.8);
	
	\draw (4,.2) --(4,-2.1);
	\draw (4,.8) --(4,1.6);
	\draw (4,2.2) --(4,2.8);
	\draw (4,3.3) --(4,3.6);
	\draw[ultra thick] (3.9,1.2) --(4.1,1.2);
	\draw[ultra thick] (3.9,2.5) --(4.1,2.5);
	
	\node (A)  at (4.7,3.1) [gray] {$\xi_1$};
	\draw[gray, ->, thin](0.2,0.2) to (2,2.2) to (4,3.1) to (A);
	\node (A)  at (4.7,3.7) [gray]  {$\xi_2$};
	\draw[gray, ->, thin](0.2,0.2) to (2,1.4) to (4,3.4) to (A);
	\node (A)  at (4.7,2.0) [gray]  {$\xi_3$};
	\draw[gray, ->, thin](0.2,0.2) to (2,1.2) to (4,2.1) to (A);
	\node (A)  at (4.7,1.6) [gray]  {$\xi_{12}$};
	\draw[gray, ->, thin](0.2,0.2) to (2,1.9) to (4,1.4) to (A);
	\node (A)  at (4.7,2.5) [gray]  {$\xi_5$};
	\draw[gray, ->, thin](0.2,0.2) to (2,2.6) to (4,2.3) to (A);
	\node (A)  at (4.7,0.2) [gray]  {$\xi_{15}$};
	\draw[gray, ->, thin](0.2,0.2) to (2,1.6) to (4,-0.3) to (A);
	\node (A)  at (4.7,-.8) [gray]  {$\xi_7$};
	\draw[gray, ->, thin](0.2,0.2) to (2,1.3) to (4,-0.0) to (A);
	
	\end{tikzpicture}
	\caption{Observations associated with the cluster point $\tilde x_2$ (solid lines) and others (dashed line). Indicated as well is the Dirichlet tessellation (Voronoi regions).\label{fig:ScenarioTree}}
\end{figure}
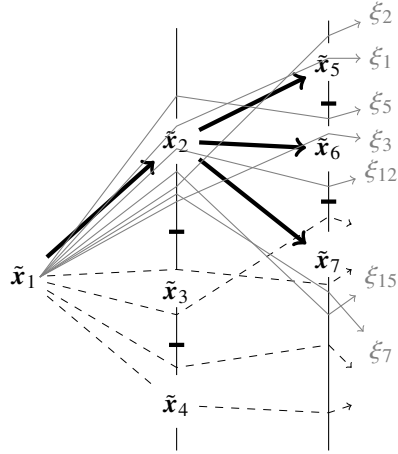

The following stage displays a difficulty, which is not present at the previous stages $t=1$ and $t=2$. Indeed, to continue the tree to stage $t=3$ it is necessary to find locations $\tilde x_i$ at the second stage and corresponding transition probabilities $\tilde p_i$. These locations and probabilities represent a \emph{conditional} distribution, conditional on the process having previously passed a particular tree state at the first stage ($\tilde x_2$, say; cf.\ \cref{fig:ScenarioTree}).

The problem encountered here is that typically \emph{none} of the observed trajectories $\xi_j$ coincides with $\tilde x_2$ exactly at the second stage, so \emph{no} data is eligible to describe the conditional distribution at the third stage. This essential difficulty is present at all subsequent stages as well (and we elaborate on this intrinsic difficulty further in \cref{sec:Trees} below).

The problem can be resolved by considering all those samples, which were associated with $\tilde x_2$ at the previous stage (cf.\ \cref{fig:ScenarioTree}). Only this subset of samples is considered further and used to build new locations and probabilities at the next stage by using the same clustering techniques as outlined for the second stage.

As a result, the approximating process has the outcomes $(\tilde x_{i_1}, \tilde x_{i_2},\tilde x_{i_3})$, where $i_2$ is a successor of $i_1$ (i.e., $i_2\in \mathcal S_{i_1}$).
This procedure can be repeated for all other locations $\tilde x_i$, $i\in \mathcal S_{i_2}$ for ${i_2}\in\mathcal S_{i_1}$ at the second stage, resulting in a tree of height $2$ and outcomes $(\tilde x_{i_1}, \tilde x_{i_2}, \tilde x_{i_3})$, where $i_1=1$ and $i_2\in \mathcal S_{i_1}$ in a cascading, nested way.
\cref{fig:ScenarioTree} displays these locations within the approximating tree including those observed scenario paths, which are associated with $\tilde x_2$. The figure displays the corresponding cluster regions as well. The clustering is called Dirichlet tessellation, or Voronoi tessellation.

\medskip
The locations at the third and later stages are found by repeatedly applying  the procedure outlined above, except that the sample paths considered further for the next $k$-means clustering are only those with common history. The paths in the approximating tree are $(\tilde x_{i_1},\dots, \tilde x_{i_T})$, where $i_{\ell+1}\in \mathcal S_{i_\ell}$ and the probabilities $\tilde p_{i_{\ell+1}}=\mathds P(\tilde X=\tilde x_{i_{\ell+1}}|\, \tilde X=\tilde x_{i_\ell})$ are conditional transition probabilities of reaching $\tilde x_{i_{\ell+1}}$ from $\tilde x_{i_\ell}$. The $k$-means problem to be solved for the node $i_t$ at stage $t$ is
\begin{equation}\label{eq:kMeans2}
\frac 1{|\mathcal N(\tilde x_{i_1},\dots,\tilde x_{i_t})|}\sum_{j \in \mathcal N(\tilde x_{i_1},\dots,\tilde x_{i_t})} \min_{\ell\in \mathcal S_{i_t}} d_{t+1}\left(\xi_{j,{t+1}},\, \tilde x_\ell\right),\qquad (\tilde x_\ell,\,\ell\in \mathcal S_{i_t}),
\end{equation}
where $d_{t+1}$ is the distance function at next stage $t+1$ and
\begin{equation*}
	\mathcal N(\tilde x_{i_1},\dots,\tilde x_{i_t}):=
	\left\{i=(i_1,\dots,i_t)\colon\ 
	\begin{aligned}
		&d_\tau(\xi_{j,\tau}, \tilde x_{i_\tau}) \le d_\tau(\xi_{j,\tau}, \tilde x_k)\\
		&\text{ for all } k\text{ with } \pred(k)=\pred(i_\tau),\ \tau\le t
	\end{aligned}
	\right\}
\end{equation*}
collects the paths associated with nearest means $(\tilde x_{i_1},\dots,\tilde x_{i_t})$ up to stage $t$. 
The total number of paths with common history is $|\mathcal N(\tilde x_{i_1},\dots,\tilde x_{i_t})|$ (see again the tie-break rule, Footnote~\ref{fn:tie-break}).

\medskip
The \emph{nested clustering} algorithm described above derives its name from this stage-wise, iterative procedure to identify an approximating process $\tilde X$ with a given branching structure.

We emphasize that clustering can be applied to data in higher dimension equally well (i.e., multivariate processes), the initial process $X$ and the approximating process $\tilde X$ do not have to be univariate (i.e., with outcomes on the real line).

Notice further that the clustering problems~\eqref{eq:kMeans} and~\eqref{eq:kMeans2} do not have to be solved to full precision in practice. Typically, a sufficient approximation is already obtained after a few subsequent iterations of the standard $k$-means clustering algorithm.
It is evident that the number $N$ of samples available should be significantly larger than $n$, the number of nodes in the approximating tree $\tilde X$. For asymptotic convergence we want to point out that the number of scenarios needed is huge, it has to hold that $N\gg n^{m_1+m_2+\dots +m_T}$, where $m_t$ is the dimension of the problem at stage $t$. This limits the use of general approximating trees in practice, especially with a higher number of stages. We refer to \citet{ShapiroNemirovski} for a more detailed analysis on the quality of approximations.

\begin{algorithm}[H] 
	\KwIn{Finite samples of $N$ trajectories, $\xi_1,\ldots,\xi_N$}
	\KwOut{Scenario tree process $\tilde X$ with a fixed branching structure}
	Stage $t=1$ is deterministic, set $\tilde X_1 = \tilde x_1$\\
	\For{$t=1,\ldots,T-1$}{
		Consider all paths $(\tilde x_{i_1},\ldots,\tilde x_{i_{t}})$ in the approximating tree
		and find new cluster means $\tilde x_\ell$ 
		at stage $t-1$ by solving the $k$-means problems
		\[
			\frac1{\left| \mathcal N(\tilde x_{i_1},\dots,\tilde x_{i_{t}})\right|}\sum_{j\in\mathcal N(\tilde x_{i_1},\dots,\tilde x_{i_{t}})} \min_{\ell \in \mathcal S_{i_{t}}} d_{t+1}\left(\xi_{j,t+1},\, \tilde x_\ell\right),\] where $d_{t+1}$ is the distance function of stage $t+1$.\\
		Define the conditional probability at stage $t$ as 
		\[\tilde p_{i_{t+1}}= \mathds P(\tilde X_{t+1}=\tilde x_{i_{t+1}}\mid \tilde x_{i_1},\dots,\tilde x_{i_{t}})=
			\frac{\left|\mathcal N(\tilde x_{i_1},\dots,\tilde x_{i_{t+1}})\right|}
			{\left| \mathcal N(\tilde x_{i_1},\dots,\tilde x_{i_{t}})\right|}.\]
	}
	\KwResult{States $\tilde x_i$ of the approximating tree $(\tilde X_1,\tilde X_2,\ldots,\tilde X_T)$ for all nodes $i\in\mathcal N$.}
	\caption{Nested clustering algorithm to determine the nodes of an approximating tree process $\tilde X$ sequentially over $T$ stages using, e.g., the $k$-means clustering algorithm.}
	\label{alg0}
\end{algorithm}

\medskip

\subsection{Stochastic approximation}\label{sec:StochApproximation}
The preceding \cref{sec:Trajectories} and the nested clustering algorithm addressed therein is based on a finite number of $N$ samples $\xi_j$, $j=1,\dots N$. In practice, a (parametric) time series model is occasionally available which allows generating arbitrarily many new sample trajectories. Every additional sample represents new information, which can and should be used to improve the approximation quality of the scenario tree. This section discusses improvements, which can be obtained by adding new samples without starting the nested clustering algorithm from scratch. The method presented here is based on \emph{stochastic approximation} (cf.\ \citet{Pflug2001} for an early application of stochastic approximation to generate scenario trees).

Suppose that an initial, approximating tree $\tilde X$ is already given. The initial tree may represent a qualified expert opinion, or may result from nested clustering. The idea outlined in this section is to modify the present tree for every  new sample path available. We use stochastic approximation to modify the actual tree given an additional, new sample.

Stochastic approximation requires us to choose a sequence of numbers $\alpha_1, \alpha_2,\dots$ beforehand, which satisfies $\alpha_k>0$, $\sum_{k=1}^\infty \alpha_k= \infty$ and $\sum_{k=1}^\infty \alpha_k^2< \infty$. A proposal, which has been proven useful in practice, is the sequence $\alpha_k=\frac 1{30+k}$, but different sequences may be more appropriate for a particular problem at hand.

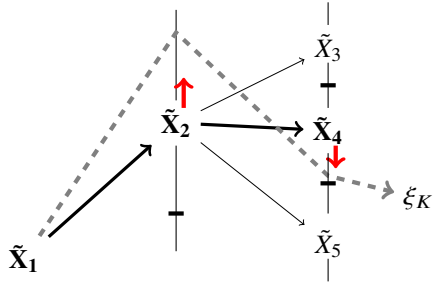
\begin{figure}[H]	
	\centering
	\begin{tikzpicture}[scale= 1.0]
	\node (1)  at (0,0.2) {$\mathbf{\tilde X_1}$};
	\node (2)  at (2,2) {$\mathbf{\tilde X_2}$};
	\node (5)  at (4,3) {$\tilde X_3$}; 
	\node (6)  at (4,1.9) {$\mathbf{\tilde X_4}$};
	\node (7)  at (4,.4) {$\tilde X_5$};
	\node (xi) at (5.2,1.0) {$\xi_{K}$};
	
	\draw[->, very thick] (1) to (2);
	\draw[->, very thin] (2) to (5);
	\draw[->, very thick] (2) to (6);
	\draw[->, very thin] (2) to (7);
	
	\draw (2,0.3) --(2,1.7);
	\draw (2,2.3) --(2,3.5);
	\draw[ultra thick] (1.9,.8) --(2.1,.8);
	
	\draw[->,gray, ultra thick, dashed](1) --(2,3.2)--(4,1.3)--(xi);
	
	\draw[->,ultra thick, red](2.1,2.2) --(2.1,2.6);
	\draw[->,ultra thick, red](4.1,1.7) --(4.1,1.4);
	
	\draw (4,.2) --(4,-0.1);
	\draw (4,.8) --(4,1.6);
	\draw (4,2.2) --(4,2.8);
	\draw (4,3.3) --(4,3.6);
	\draw[ultra thick] (3.9,1.2) --(4.1,1.2);
	\draw[ultra thick] (3.9,2.5) --(4.1,2.5);
	\end{tikzpicture}
	\caption{Every stochastic approximation step modifies one path within the tree (thick) slightly towards the new observation $\xi_{K}$ (dashed), cf.~\eqref{eq:convex}.\label{fig:StochasticApproximation}} 
\end{figure}

Suppose an additional trajectory is available, which we denote by $\xi_{k}=(\xi_{k,1},\dots,\xi_{k,T})$ for convenience. We modify the given tree $\tilde X^{(k-1)}$ based on the new information $\xi_{k}$ to get a new stochastic tree $\tilde{X}^{(k)}$. This is achieved by identifying the branch in the tree $\tilde X^{(k)}(i_1)$ which is closest to $\xi_{k,1}$ at the 1\textsuperscript{st} stage, then the node $i_2$ in the subtree which is closest to $\xi_{k,2}$ at the 2\textsuperscript{nd} stage, etc. Suppose the path $\left(\tilde X^{(k-1)}(i_1),\dots, \tilde X^{(k-1)}(i_T)\right)$  is sequentially closest to the observation $\xi_k$. The stochastic approximation step then updates the tree $\tilde X^{(k-1)}$ by shifting this path slightly. Each individual component of the trajectory is updated to
\begin{equation}\label{eq:convex}
\tilde X^{(k)}(i) \leftarrow (1-\alpha_{k})\cdot\tilde X^{(k-1)}(i)+ \alpha_{k}\cdot\xi_{k,i}
\end{equation}
for every $i\in \{i_1,\dots, i_T\}$.
Note that~\eqref{eq:convex} is a convex combination, shifting the path $\tilde X^{(k)}(i_t)$ towards the new observation $\xi_k$ with a force specified by $\alpha_{k}$.

\cref{fig:StochasticApproximation} illustrates a modifying step of the stochastic approximation procedure, the path is shifted along the horizontal arrows in direction of the new sample $\xi_{k}$.

The stochastic approximation procedure repeats the modifying step with new scenarios successively, until satisfactory convergence is obtained. We refer to \citet{Kushner2003} for details on stopping criteria.

\medskip
\citet{PflugPiSzenarioGen} present an extensive explanation for an algorithm for generating scenario trees. They also employ the concept of stochastic approximation procedure discussed in \cref{sec:StochApproximation}. The following example demonstrates generation of scenario trees using the concept presented therein.

\begin{example}[Running maximum]
	Consider the running maximum process in \cref{fig:rmp} in 4 stages. This process is given by 
	\begin{equation}\label{eq:rmp}
	M_t = \max \Big\{ \sum_{i=1}^{t^\prime} \xi_i : t^\prime \leq t\Big\},\ t=1,\dots,4\ \text{ with} \ \xi_i \sim \mathcal N(0,1) \ \text{i.i.d.}
	\end{equation} 
	The process $M_t$ is non-Markovian since it depends on the history. Therefore, we can use this process to illustrate  generation of scenario trees with a fixed branching structure and a certain number of iterations for the stochastic approximation algorithm.
	\begin{figure}[H]
		\centering
		\includegraphics[height=6cm,width=10cm]{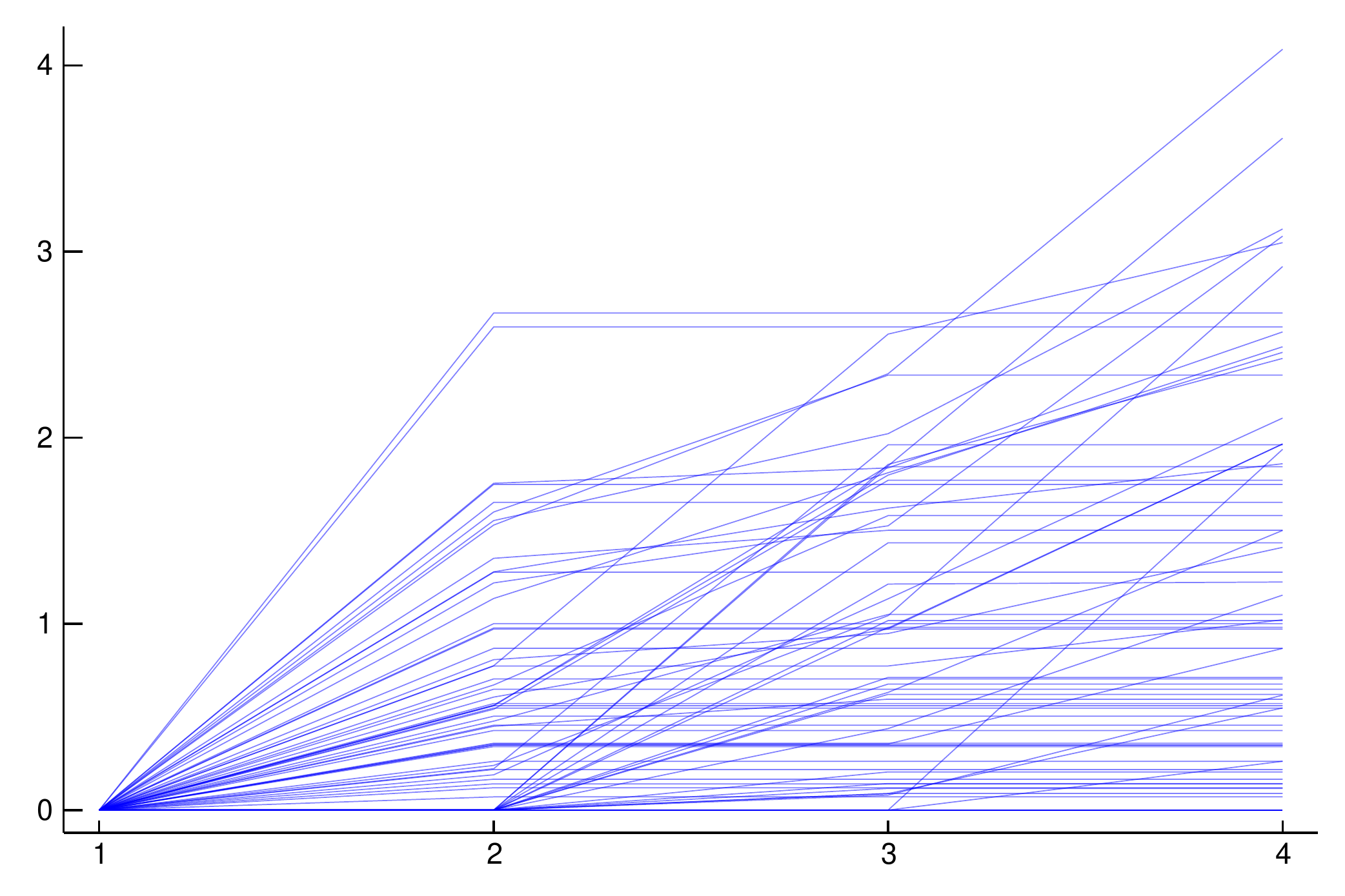}
		\caption{Scenarios of a running maximum process in four stages}
		\label{fig:rmp}
	\end{figure}
\cref{fig:rmp2} below shows two scenario trees with branching structure\footnote{We refer to the number of nodes at each stage, $(|\mathcal N_1|,|\mathcal N_2|,\dots, |\mathcal N_T|)$, as the \emph{branching structure} of the process.} $(1,2,2,2)$ and $(1,3,3,3)$ approximating this process.
They are generated with \texttt{tree\_approximation()} function contained in the package
\href{https://kirui93.github.io/ScenTrees.jl/stable/}{\texttt{ScenTrees.jl}}.\textsuperscript{\ref{fn:ScenTrees}}
	\begin{figure}[H]
		\centering
		\begin{subfigure}[t]{0.49\textwidth}
			\includegraphics[width=\textwidth, trim=0 0 0 20, clip]{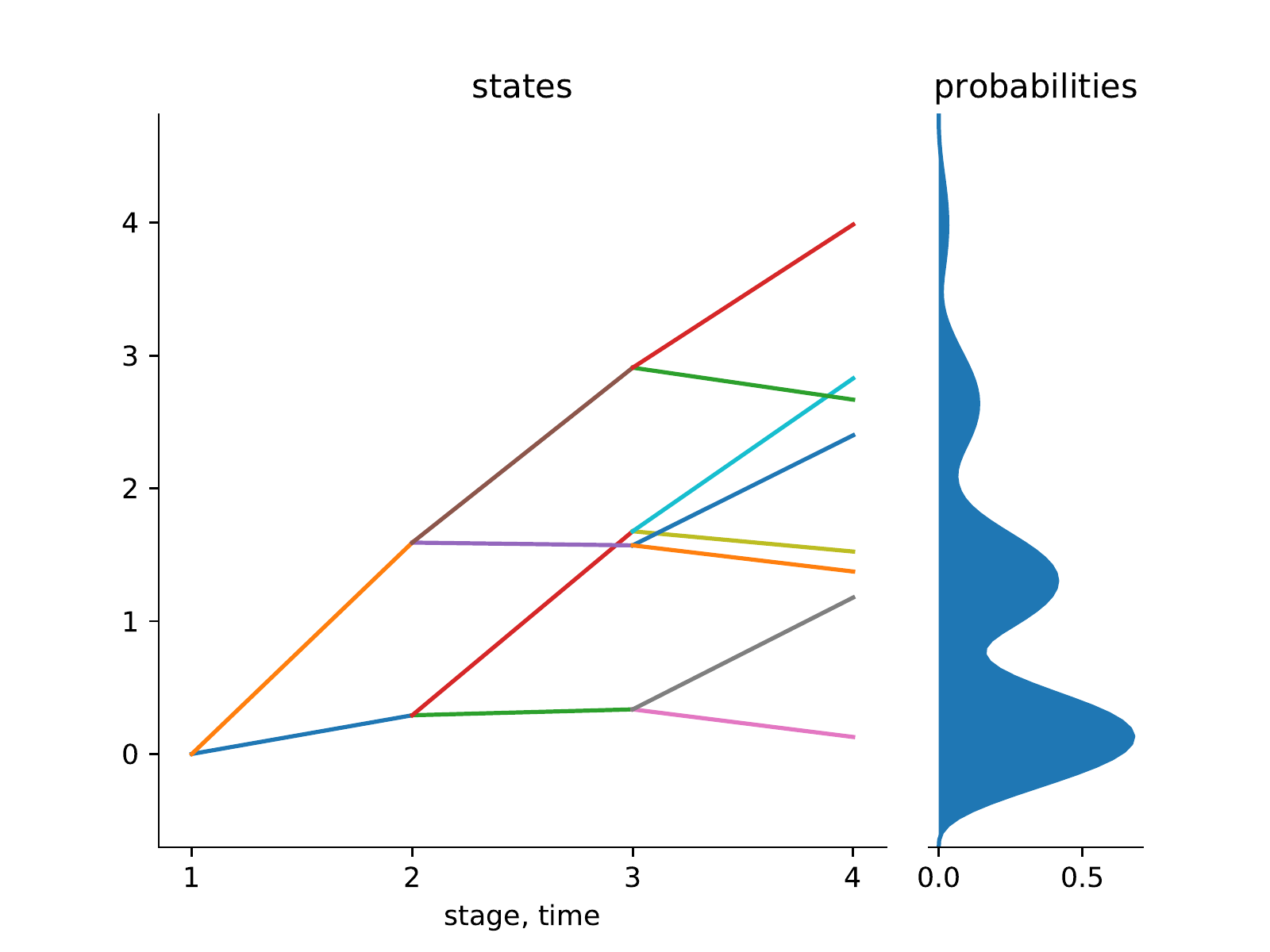}
			\caption{Binary tree approximating process in \cref{fig:rmp} with a transportation distance of $0.313$}
			\label{fig:rmtree2}
		\end{subfigure}\hfill
		\begin{subfigure}[t]{0.49\textwidth}
			\includegraphics[width=\textwidth, trim=00 00 00 20, clip]{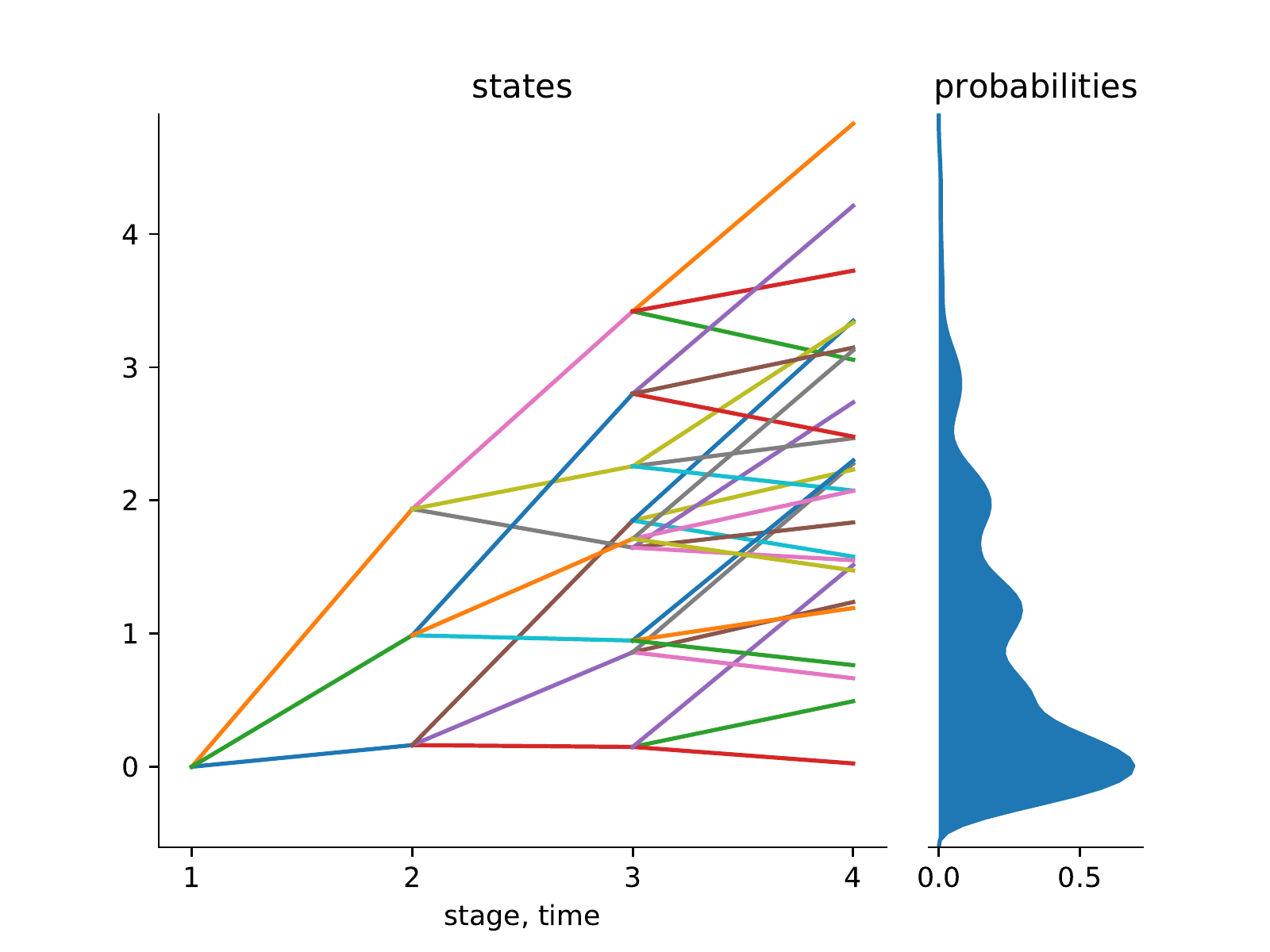}
			\caption{Tree with 3 branches at each node approximating process in \cref{fig:rmp} with a transportation distance of $0.107$}
			\label{fig:rmtree3}
		\end{subfigure}
		\caption{Example of two scenario trees with different branching structures approximating process in \cref{fig:rmp}.}
		\label{fig:rmp2}
	\end{figure}
	\cref{fig:rmtree3} has a better approximating quality than \cref{fig:rmtree2}. This is because \cref{fig:rmtree3} has more branches than \cref{fig:rmtree2} and therefore represents much more information. This can also be seen through the difference in the multistage distance. The multistage distance between the origin maximum running process $M_t$ (with realization in \cref{fig:rmp}) and \cref{fig:rmtree3} is $0.107$ while that of \cref{fig:rmtree2} is $0.313$.
\end{example}

\section{Scenario tree generation with limited data}\label{sec:Trees}
In real world applications it is not unusual that data available is limited, i.e., that the number of available samples $N$ is fixed and not large; \citet{Woodruff2015} report further on this difficulty. To generate a stochastic tree based on limited data it is necessary to modify the approach in comparison to the \cref{sec:Trajectories} above. In these previous sections we combine and compress scenarios with similar outcomes at each step to locate new, continuative nodes of the tree. These methods need many, or unlimited data to specify a good approximating tree, and thus cannot be applied in case of limited data.

In case of limited data it is necessary to learn as much as possible from the observations available. The methods presented below differ in the problem already mentioned in \cref{sec:Trajectories}, that is, how to deal with a specific history, which is not exactly met by any data observed.

To overcome this intrinsic difficulty we generate new and additional, but different samples based on the fixed number of observed samples $\xi_j$, $j=1,\dots N$, available. This is accomplished by estimating the density describing the transition.

\medskip
These additional samples can then be used to apply the procedures outlined above, but in what follows we describe further algorithms which make direct use of the estimated transition distributions as well.

\subsection{Scenario generation from limited data}
In contrast to nested clustering (cf.\ \cref{sec:Trajectories})
we estimate the distribution of the transition given the \emph{precise} history. This distribution can be estimated by non-parametric kernel density estimation. That is, the density at stage $t+1$ conditional on the history $x_{[t]}:=(x_1,x_{2},\dots,x_{t})$ is estimated by 
\begin{equation}\label{eq:Conditioal}
	\hat f_{t+1}(x_{t+1}\mid x_{[t]})=\sum_{j=1}^N  w_j(x_{[t]}) \cdot k_{h_N}\Big(x_{t+1}-\xi_{j,t+1}\Big),
\end{equation}
where $k(\cdot)$ is a kernel function, $h_N$ is the bandwidth and $k_{h_N}(\cdot):=\frac{1}{h_N^{m_{t+1}}} k\left(\frac{\cdot}{h_N}\right)$. The weights
\begin{equation}\label{eq:weight}
w_j(x_{[t]})=\frac{k\left(\frac{x_{i_1}-\xi_{j,1}}{h_N}\right)\cdot\ldots\cdot k\left(\frac{x_{i_t}-\xi_{j,t}}{h_N}\right)}{\sum_{\ell=1}^N k\left(\frac{x_{i_1}-\xi_{\ell,1}}{h_N}\right)\cdot\ldots\cdot k\left(\frac{x_{i_t}-\xi_{\ell,t}}{h_N}\right)}
\end{equation}
depend on the history $x_{[t]}$ (for consistency, we set $w_j(x_{[1]}):=\frac 1 N$ for $t=1$).

It can be shown that the approximation is asymptotically optimal for \\ $h_N\sim N^{-1/(m_1+\dots+m_t+4)}$ (cf.\ Silverman's rule of thumb, see also \citet{PflugPichler2016}). Note as well that the weights sum to $1$ in~\eqref{eq:weight} ,i.e., for every $x_{[t]}$, $\sum_{j=1}^N w_j(x_{[t]})=1$. Further note that the estimator~\eqref{eq:Conditioal} for the conditional density involves \emph{all} samples $\xi_j$, $j=1,\dots, N$.
The weight $w_j(x_{[t]})$ is further large, if the first $t$ stages $(\xi_{j,1},\dots, \xi_{j,t})$ of the observation $\xi_j$ are similar to $x_{[t]}$, and negligible otherwise. In this way the weight $w_j(x_{[t]})$ determines the importance of the corresponding observation $\xi_j$ for $x_{[t]}$, and the estimator~\eqref{eq:Conditioal} takes this information adequately into account.

Univariate Kernel functions, which have proven useful are the Epanechnikov kernel $k(x)= \max\left\{\frac 3 4(1-x^2),\,0\right\}$ and the logistic kernel $k(x)= \frac 2{(e^x+e^{-x})^2}$, cf.\ \citet{Tsybakov}. For multivariate data (i.e.\ higher dimension, $m_t>1$) univariate kernels can be multiplied to obtain a multivariate kernel function, $k(x_1,\dots, x_{m_t})= k(x_1)\cdot\ldots\cdot k(x_{m_t})$.

\paragraph{Samples from the conditional distribution~\eqref{eq:Conditioal}~--- the composition method.}
The composition method allows to quickly find samples from the conditional distribution $f_{t+1}(\,\cdot\,|\,x_{[t]})$ in~\eqref{eq:Conditioal}. Indeed, pick a random number $U\in (0,1)$, where $U$ is uniformly distributed. Then there is a summation index $j^*\in \{1,2,\dots N\}$ so that
\begin{equation}\label{compeq}
	\sum_{j=1}^{j^*-1}w_j(x_{[t]})< U\le \sum_{j=1}^{j^*}w_j(x_{[t]}).
\end{equation}
Note, that $j^*$ is distributed on $\{1,2,\dots N\}$ with probability mass function $P(j^*= j)= w_j(x_{[t]})$.
For this index $j^*$ then randomly pick an instant from the kernel density $\frac{1}{h_N^{m_{t+1}}} k\left(\frac{\,\cdot\ \, - \xi_{j^*,t+1}}{h_N}\right)$ shifted to $\xi_{j^*,\,t+1}$, cf.~\eqref{eq:Conditioal}. This sample is finally distributed with the desired density $\hat f_{t+1}(\,\cdot\,|\, x_{[t]})$. The composition method thus is a simple, quick and cheap method to make samples from the additively composed density~\eqref{eq:Conditioal} accessible.

\subsection{A direct tree generation approach with limited data}
Note that~\eqref{eq:Conditioal} describes a continuous distribution. To specify the tree process we may approximate the distribution again by finding a discrete distribution with locations and corresponding weights, which approximate~\eqref{eq:Conditioal} adequately.

A scenario tree can be generated directly by exploiting the conditional density estimates introduced above. To this end approximate the distribution in~\eqref{eq:Conditioal} with density $\hat f_1(\cdot)$ by the discrete distribution $\sum_{\ell\in \mathcal S_1} \tilde p_\ell \cdot \delta_{\tilde x_\ell}$ located at $\tilde x_\ell$ ($\ell\in \mathcal S_1$), where each $\tilde x_\ell$ is reached with probability $\tilde p_\ell$. These locations constitute the states at stage $t=1$, building the origin of the approximating tree.

For the next step define $\mathbf{\tilde{x}}_{i_1}:=(\tilde x_{i_1},\tilde x_{i_2})$ and consider the distribution with density $f_1(\cdot|\,\mathbf{\tilde x}_{i_1})$. By applying the same approximation procedure again one finds a discrete, approximating distribution $\sum_{\ell\in \mathcal S_{i_1}} p_\ell\cdot \delta_{\tilde x_\ell}$. The tree then is continued with $(\tilde x_{i_1}, \tilde x_{i_2},\tilde x_{i_3})$, where $i_2\in \mathcal S_{i_1}$ is a successor of $i_1$.

A tree is finally found by repeatedly applying the method outlined above at every node in the subtree already available.
\medskip

\begin{algorithm}[H]
	\SetKwData{Init}{\textbf{Initialization:}}
	\SetKwData{Iter}{\textbf{Iteration:}}
	\KwIn{Samples $\xi_j$ with $j=1,\ldots,N$ and a kernel function $k(\cdot)$}
	\KwOut{A random vector $x=(x_1,\ldots,x_T)$ following the distribution in~\eqref{eq:Conditioal}}
	\Init $w=(1.0,\ldots,1.0)$ is an $N$-dimensional vector of sample weights, cf.~\eqref{eq:weight}.\\
	\For{$t=1,\ldots,T$}{
		Normalize the weights, set $w_j = \frac{w_j}{\sum_{\ell=1}^{N} w_\ell}$.\\
 		Set the bandwidth $h_t = \sigma_t \cdot N_t ^{-\frac{1}{m_t+4}}$, where $m_t$ is the dimension, $N_t = \frac{\big(\sum_{j=1}^{N} w_j\big)^2}{\sum_{j=1}^{N} w_{j}^{2}}$ is the effective sample size and $\sigma_{t}^{2} = \var\big(\xi_{j,t}\colon j=1,\ldots,N\big)$ is the variance at stage $t$.
 		\\
		Draw a uniform random number $U \in [0,1]$.\\
		Find an index $j^\ast$ such that $\sum_{j=1}^{j^\ast-1} w_j < U \leq \sum_{j=1}^{j^\ast} w_j$ (cf.~\eqref{compeq}).\\
 		Pick a random instance $K_t$ from the distribution with density $k(\cdot)$ and set $x_t = \xi_{j^\ast,t} + h_t \cdot K_t$.\\
		Re-calculate the sample weights for the next stage:\\
		\eIf{Markovian}{
			\For{$j=1,\ldots,N$}{$w_j = k_{h_t}(x_t - \xi_{j,t})$}}{
			\For{$j=1,\ldots,N$ }{$w_j = w_j \cdot k_{h_t}(x_t - \xi_{j,t})$}}}
		\KwResult{A random trajectory $(x_1,\dots,x_T)$ following the distribution~\eqref{eq:Conditioal}}
	\caption{Generate a new trajectory using conditional density estimation. 
	\href{https://kirui93.github.io/ScenTrees.jl/stable/}{\texttt{ScenTrees.jl}\textsuperscript{\ref{fn:ScenTrees}}} provides this function as \texttt{kernel\_scenarios().}
	}
	\label{alg3}
\end{algorithm}

\subsection{Stochastic approximation with limited data}
The stochastic approximation algorithm outlined in \cref{sec:StochApproximation} bases on additional sample paths, and each new sample path is used to improve the scenario tree. To obtain convergence of the method it is necessary to have an unlimited number of sample paths available, but this is certainly not the case for a limited number of data. However, the conditional density estimation methods described in this section allow generating new sample paths (cf.\ also \citet{Haerdle1997}).

\paragraph{New sample paths from observations.}
Every new sample path starts with $x_1$ at the first stage. Using~\eqref{eq:Conditioal} and the composition method one may find a new sample $x_2$ from $\hat f_2(\cdot\mid x_[1])$ at the second stage. Then generate another sample $x_3$ from $\hat f_3(\cdot\mid x_{[2]})$ at the third stage, a new sample $x_t$ from $\hat f_t(\cdot\mid x_{[t-1]})$ at the stage $t$, etc. Iterating the procedure until the final stage $T$ reveals a new sample path $x=(x_1,x_2,\dots,x_T)$, generated from the initial data $\xi_1,\dots,\xi_N$ directly~(cf.~\cref{alg3}).

\medskip
Generating new scenario paths can be repeated arbitrarily often to get new, distinct scenario paths. Each new sample path $x$ can be handed to the stochastic approximation algorithm (\cref{sec:StochApproximation}) to improve iteratively a tentative, approximating tree.

With this modification the stochastic approximation algorithm in \cref{sec:StochApproximation} is eligible even in case of limited data.

\begin{example}[Running Maximum using~\cref{alg3}]
	Consider a case where we have just $N=100$ scenarios of running maximum data. We use this limited data and learn as much as possible so that we can generate new and additional samples based on this data. These generated samples are fed directly into stochastic approximation procedure to generate a scenario tree. 
	\cref{fig:rmdg} shows 100 trajectories generated from the data using \cref{alg3} and the resulting scenario tree with a branching structure of $(1,3,3,3)$.
	\begin{figure}[H]
		\centering
		\begin{subfigure}[t]{0.49\textwidth}
			\includegraphics[width=\textwidth, trim=0 0 0 20, clip]{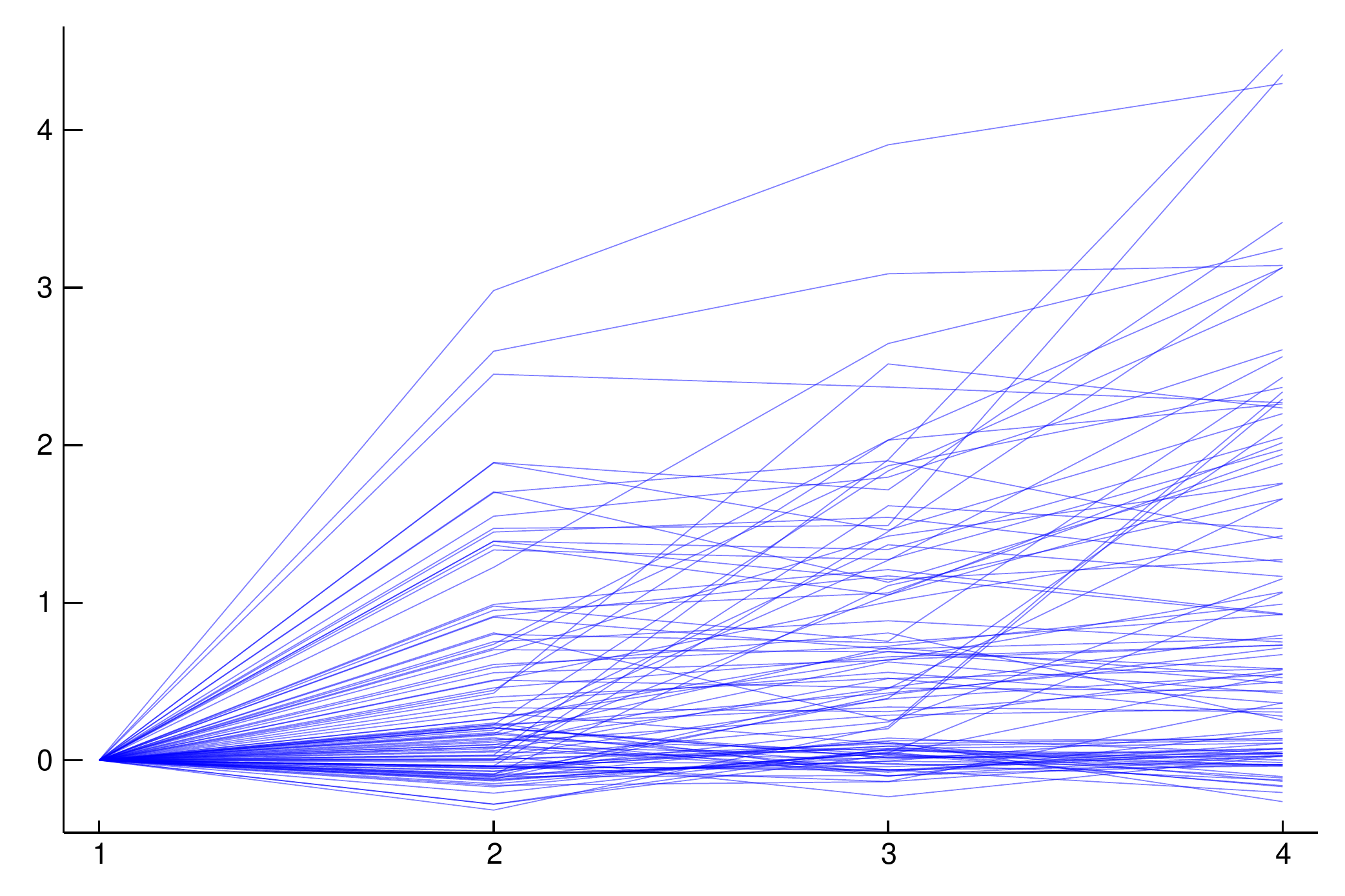}
			\caption{100 trajectories generated from the process in \cref{fig:rmp} using \cref{alg3}}
			\label{fig:rmdatag}
		\end{subfigure}\hfill
		\begin{subfigure}[t]{0.49\textwidth}
			\includegraphics[width=\textwidth, trim=00 00 00 20, clip]{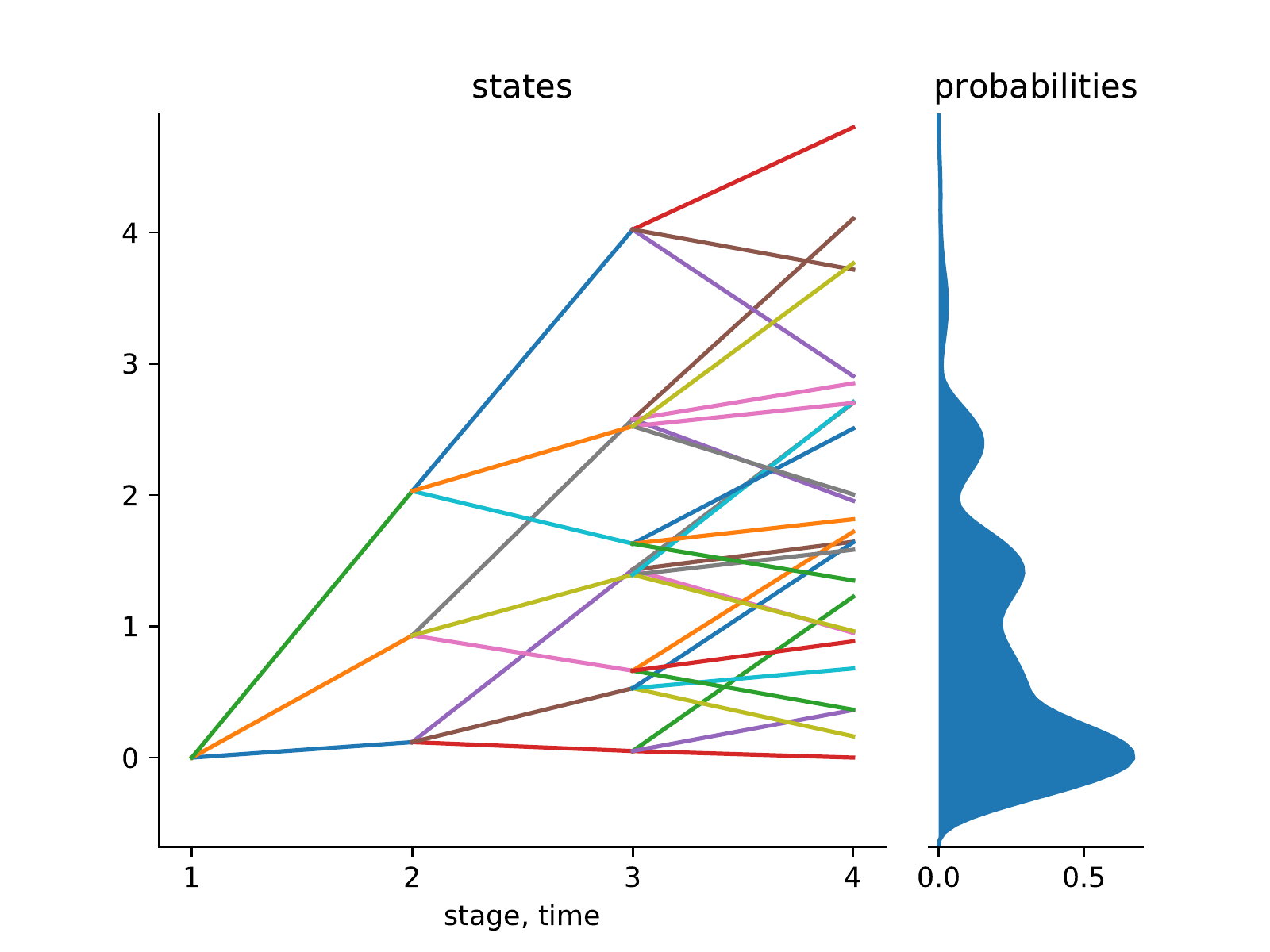}
			\caption{Tree with 3 branches at each node approximating process in \cref{fig:rmdatag} with distance of $0.114$}
			\label{fig:rmtreeg}
		\end{subfigure}
		\caption{100 trajectories generated from a limited running maximum data and a resulting scenario tree, generated with
		\href{https://kirui93.github.io/ScenTrees.jl/stable/}{\texttt{ScenTrees.jl}\textsuperscript{\ref{fn:ScenTrees}}}}
		\label{fig:rmdg} 
	\end{figure}
	It is clear that the data generated (cf.\ \cref{fig:rmdg}) is similar to the original data in \cref{fig:rmp}. Therefore, this method proves efficient to generate new data for generation of scenario trees using stochastic approximation procedure (cf.\ \cref{fig:rmtreeg}).
\end{example}

\section{Generation of scenario lattices} \label{sec:Lattice}
In many situations of practical relevance the process to be approximated by a tree is known to have special properties. Of course, this information should be exploited when constructing the approximating tree in order to reduce the computational effort and time, and to obtain trees with better approximation quality.

\medskip
In what follows we discuss the situations of Markovian processes and processes with fixed states separately. A main characteristic of trees in finance is arbitrage. For this we discuss arbitrage opportunities as well, and we provide a new result for trees with stochastic dominance in addition.

\subsection{Markovian processes and lattices}\label{sec:Markov}

Estimating a lattice follows the same principle as estimating a tree. However, the lattice generation procedure is of significantly lower complexity. If $D_t =\{\tilde{x}_{t,1}, \dots, \tilde{x}_{t,i_t} \}$ is the set of states at stage $t$, then the induced distribution between the conditional lattice given $\tilde{x}_{t,i}$ and the Markovian process $\xi_t$ must be calculated only for all nodes from $D_t$ and not for all paths leading to them.
As in the tree case, the calculation of the transport distance  requires a backward algorithm but can be based on the approximation of the one-period transition probabilities  $P(\xi_{t+1}\in A \mid \xi_t = x)$.
The average aberration distance\footnote{This aberration distance incorporates the additional parameter $r\ge 1$.}
\begin{equation}\label{aad}
	\left(\E\left(\sum_{t=1}^T \,d_t\big(\xi_t, \mathbf T_t(\xi_t)\big)\right)^r\right)^{\nicefrac1r}
\end{equation}
(cf.~\eqref{eq:distT}) is an upper bound for the nested distance, but this upper bound may be conservative.

One has to consider the full conditional future process for the estimation of a scenario lattice, however, one does not have to consider the past. This fact allows the following
simplifications for the algorithms:
\begin{enumerate}[label=(\roman*)]
	\item When applying the nested clustering or stochastic approximation (\cref{sec:Trajectories,sec:StochApproximation} above), for example, all samples can be reused at intermediate stages, as they do not have to have the same history.
	\item When applying the density estimator~\eqref{eq:Conditioal}, the weights in~\eqref{eq:weight} simplify to
	\begin{equation*}
	w_j(x_{[t]})=\frac{k\left(\frac{x_{i_t}-\xi_{j,t}}{h_N}\right)}{\sum_{\ell=1}^N k\left(\frac{x_{i_t}-\xi_{\ell,t}}{h_N}\right)}.
	\end{equation*}
\end{enumerate}

Applications of lattice scenario  models  in an hydro storage management system can be found, e.g., in \citet{Wozabal2013}.
\medskip
\begin{example}[Gaussian random walk]
Consider the Gaussian random walk process $X_n = \sum_{k=1}^{n} Y_k$  on 5 stages, where the random variables $\big(Y_k\big)_{k\geq1}$ are i.i.d.\@ and the distribution of each $Y_k$ is normal.
Simple examples of this process would be a path traced by a molecule as it travels through a liquid or a gas, the price of a fluctuating stock and the financial status of a gambler.
Random walks are also simple representations of the Markov processes.
\cref{fig:gaussianL} employs \cref{alg2} to approximate this process using a scenario lattice.
\begin{figure}[H]
	\centering
	\begin{subfigure}[t]{0.49\textwidth}
		\includegraphics[width=\textwidth, trim=0 0 0 20, clip]{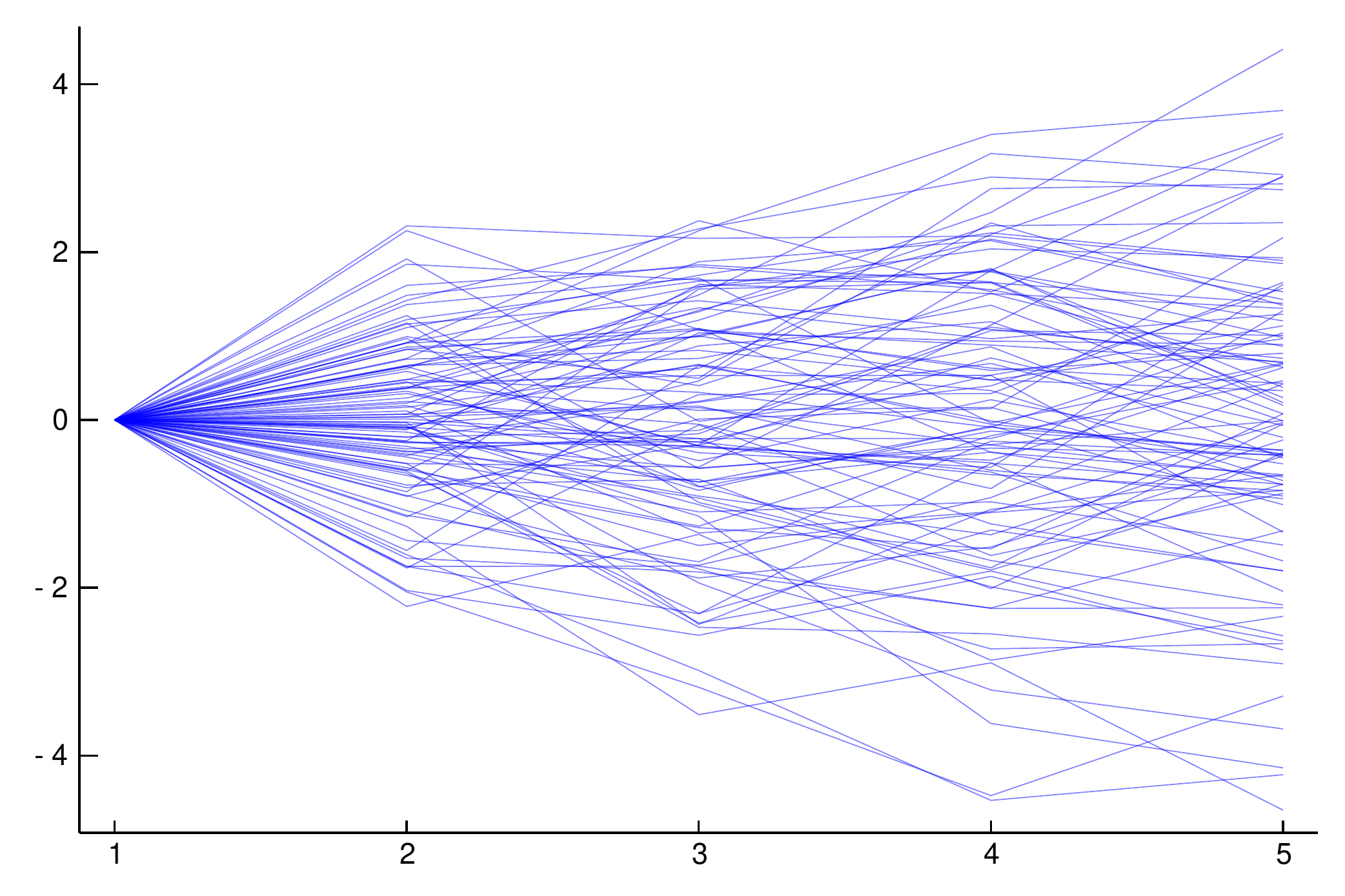}
		\caption{100 samples from Gaussian random walk process}
		\label{fig:gaussian}
	\end{subfigure}\hfill
	\begin{subfigure}[t]{0.49\textwidth}
		\includegraphics[width=\textwidth, trim=00 00 00 20, clip]{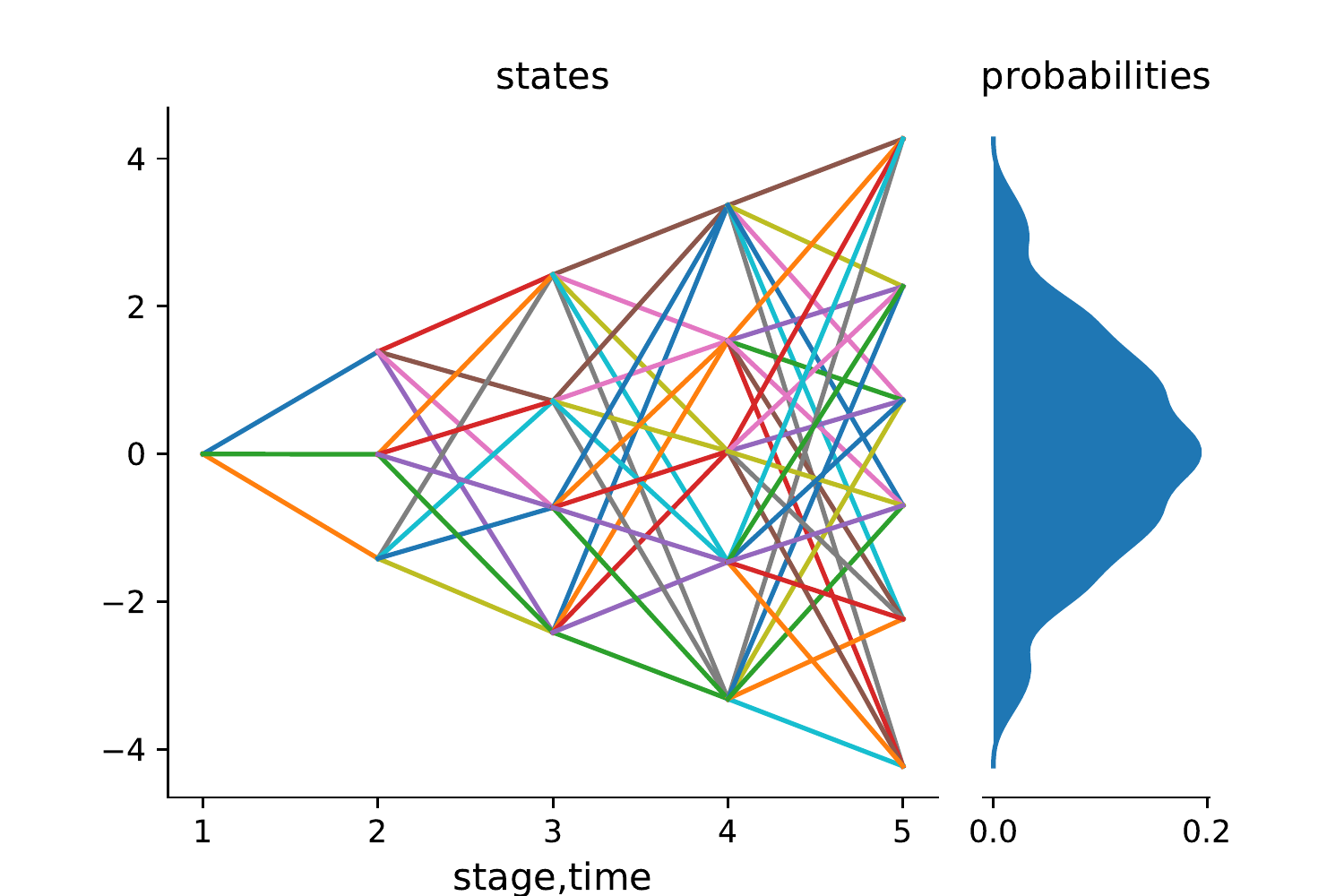}
		\caption{Scenario lattice with branching structure $(1,3,4,5,6)$ approximating process in \cref{fig:gaussian} with a distance of $0.554$} 
		\label{poissonLattice}
	\end{subfigure}
	\caption{100 samples from Gaussian random walk process and the resultant scenario lattice, generated with \href{https://kirui93.github.io/ScenTrees.jl/stable/}{\texttt{ScenTrees.jl}}}
	\label{fig:gaussianL}
\end{figure}
\end{example}

\begin{algorithm}[H] 
	\SetKwData{Init}{\textbf{Initialization:}}
	\SetKwData{Iter}{\textbf{Iteration:}}
	\KwIn{Let $K$ be the number of iterations.}
	\KwOut{Scenario lattice with values $\tilde X^{(K)}$ and transition probabilities.}
	\Init Set the transportation distance $c_E=0$, set the counters $c(n) = 0$ for all nodes and let $\tilde X^{(0)}$ be an initial lattice chosen by expert opinion, e.g. Also, choose a non-negative and non-increasing sequence $\alpha_{k}$ such that $\sum_{k=1}^{\infty} \alpha_{k} = \infty$ and $\sum_{k=1}^{\infty} \alpha_{k}^{2} < \infty$.\\
	\For{$k=1,\ldots,K$}{
	 Use a new and independent trajectory $\xi^{(k)} = (\xi_{1}^{(k)},\ldots,\xi_{T}^{(k)})$\\
	Find a scenario  $s=(i_1,\ldots,i_T)$ of nodes with the closest lattice entry such that 
	\begin{equation}\label{eq:si}
		i_t\in \argmin_{n^\prime \in \mathcal{N}(t)} d_t\big(\tilde X^{(k)}(n^\prime),\xi_{t}^{(k)}\big),
	\end{equation} where $\mathcal{N}(t)$ are all nodes in stage $t$.\\
	\For{$t =1,\ldots,T$}{
	Increase all counters $c(i_t) = c(i_t) + 1$ for all nodes in $s$ (cf.~\eqref{eq:si}).\\
	Modify the values of the nodes in these entries to be \[\tilde X^{(k)} (i_t) = \tilde X^{(k-1)}(i_t) - \alpha_{k} \cdot r\cdot d_t\big(\tilde X^{(k-1)}(i_t),\xi_{t}^{(k)}\big)^{r-1}\cdot \big(\tilde X^{(k-1)}(i_t) - \xi_{t}^{(k)}\big).\]
	The values on other node entries remain unchanged.
}
	Set $c_E = c_E + \big(\sum_{t=1}^{T}d_t(\tilde X^{(k-1)}(i_t),\xi_{t}^{(k)})\big)^r.$}
	\KwResult{Set the conditional probabilities $p(n) = \frac{c(n)}{K}$. The quantity approximating $\big(\E d(\xi,T(\xi))^r\big)^{\nicefrac1r}$ (cf.~\eqref{aad}) is given by $\big(c_E/ K\big)^{\nicefrac1r}$.}
	\caption{Generating a scenario lattice with a fixed branching structure using stochastic approximation.
	\href{https://kirui93.github.io/ScenTrees.jl/stable/}{\texttt{ScenTrees.jl}\textsuperscript{\ref{fn:ScenTrees}}} provides this algorithm as \texttt{lattice\_approximation()}}
	\label{alg2} 	
\end{algorithm}

\subsection{Processes with fixed states}\label{sec:FixedStates}
In selected situations of practical relevance the distribution at stage $t$ is known approximately, or the support of this distribution or, for example, its quantiles can be described  beforehand. In these situations it might make sense to fix the states beforehand and to use the algorithms introduced above to estimate the transition probabilities only, without modifying the states. Examples of such processes include mean reverting processes.

\begin{figure}[H]
	\centering
	\begin{tikzpicture}[scale= 0.85]
	\node (1)  at (0,0) {$\tilde x_1$};

	\node (2)  at (2,1.7) {$\tilde x_2$};
	\node (3)  at (2,1) {$\tilde x_3$};
	\node (4)  at (2,0) {$\tilde x_4$};
	\node (5)  at (2,-1){$\tilde x_5$};
	\node (6)  at (2,-1.7){$\tilde x_6$};

	\node (11)  at (4,1.7) {$\tilde x_7$};
	\node (10)  at (4,1) {$\tilde x_8$};
	\node (9)  at (4,0) {$\tilde x_9$};
	\node (8) at (4,-1){$\tilde x_{10}$};
	\node (7) at (4,-1.7){$\tilde x_{11}$};

	\foreach \x in {1.6,1,0,-1,-1.6}\node at (5,\x) {\dots};

	\node (16) at (6,1.7) {$\tilde x_.$};
	\node (15) at (6,1) {$\tilde x_.$};
	\node (14) at (6,0) {$\tilde x_.$};
	\node (13) at (6,-1){$\tilde x_.$};
	\node (12) at (6,-1.7){$\tilde x_.$};

	\node (21) at (8,1.7) {$\tilde x_.$};
	\node (20) at (8,1) {$\tilde x_.$};
	\node (19) at (8,0) {$\tilde x_.$};
	\node (18) at (8,-1){$\tilde x_.$};
	\node (17) at (8,-1.7){$\tilde x_.$};

	\foreach \x in {2,...,6}                            \draw[->, ultra thin] (1) to (\x);
	\foreach \x in {2,...,6}  \foreach \y in {7,...,11} \draw[->, ultra thin] (\x) to (\y);
	\foreach \x in {12,...,16}\foreach \y in {17,...,21}\draw[->, ultra thin] (\x) to (\y);

	\node (t0) at (0.3,-2.5) {$t=1$};
	\node (t1) at (2,-2.5) {$t=2$};
	\node (t2) at (4,-2.5) {$t=3$};
	\node (t3) at (7.7,-2.5) {$t=T$};
	\end{tikzpicture}
	\caption{A scenario lattice process with fixed states\label{fig:Fixed}}
\end{figure}
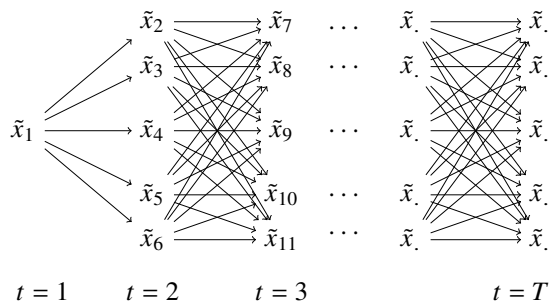
\Cref{fig:Fixed} displays a scenario lattice with 5, fixed states, which do not change over time. The points $\tilde x_2$ to $\tilde x_6$ can be determined beforehand, and they are kept fixed over time. For a univariate distribution they can be chosen to be the quantiles, e.g. That is,
$x_2$ is the $\nicefrac 1{10}$-quantile,
$x_3$ the $\nicefrac 3{10}$-quantile, \dots\ and
$x_6$ the $\nicefrac 9{10}$-quantile, resp.

\section{Trees lacking a predefined structure}\label{sec:Treeswithout}
Up to now we have assumed that the structure of the tree (i.e., the number of nodes and the predecessor vector $\pred$) has been chosen and is fixed. This a quite restrictive, since hardly anybody can decide beforehand what tree topology bes suits a given decision problem. Here is how to generate a tree with a given approximation quality without choosing the tree topology beforehand.
\begin{enumerate}[label=(\roman*)]
\item Choose (for each stage) a maximally acceptable distance between the estimated conditional densities and the discrete approximations.
\item Start with a simple tree (e.g., a binary one).
\item If during the tree construction the required distance cannot be achieved with the given number of successors, increase the branching factor by 1 until the distance condition can be satisfied.
\end{enumerate}
Details of this procedure can be found in \citet{PflugPichlerBuch}.
\section{Application of scenario lattice generation methods to an electricity pricing data}\label{sec:application}

To summarize the methods discussed in the preceding sections, we consider a limited electricity load data recorded for each hour of the day for the year 2017. The goal in what follows is to approximate the load of electricity for each hour of the day from this data using a scenario lattice. The scenario lattice has the same number of stages as the number of hours in the data. It also has a certain fixed branching structure (i.e., the number of nodes connected to each succeeding nodes in the lattice).

This computation has been done using the
\href{https://kirui93.github.io/ScenTrees.jl/stable/}{\texttt{ScenTrees.jl}}\textsuperscript{\ref{fn:ScenTrees}}  
package which has been tested for Julia~$\geq 1.0$.

\begin{figure}[H]
	\centering
	\includegraphics[width=0.75\textwidth]{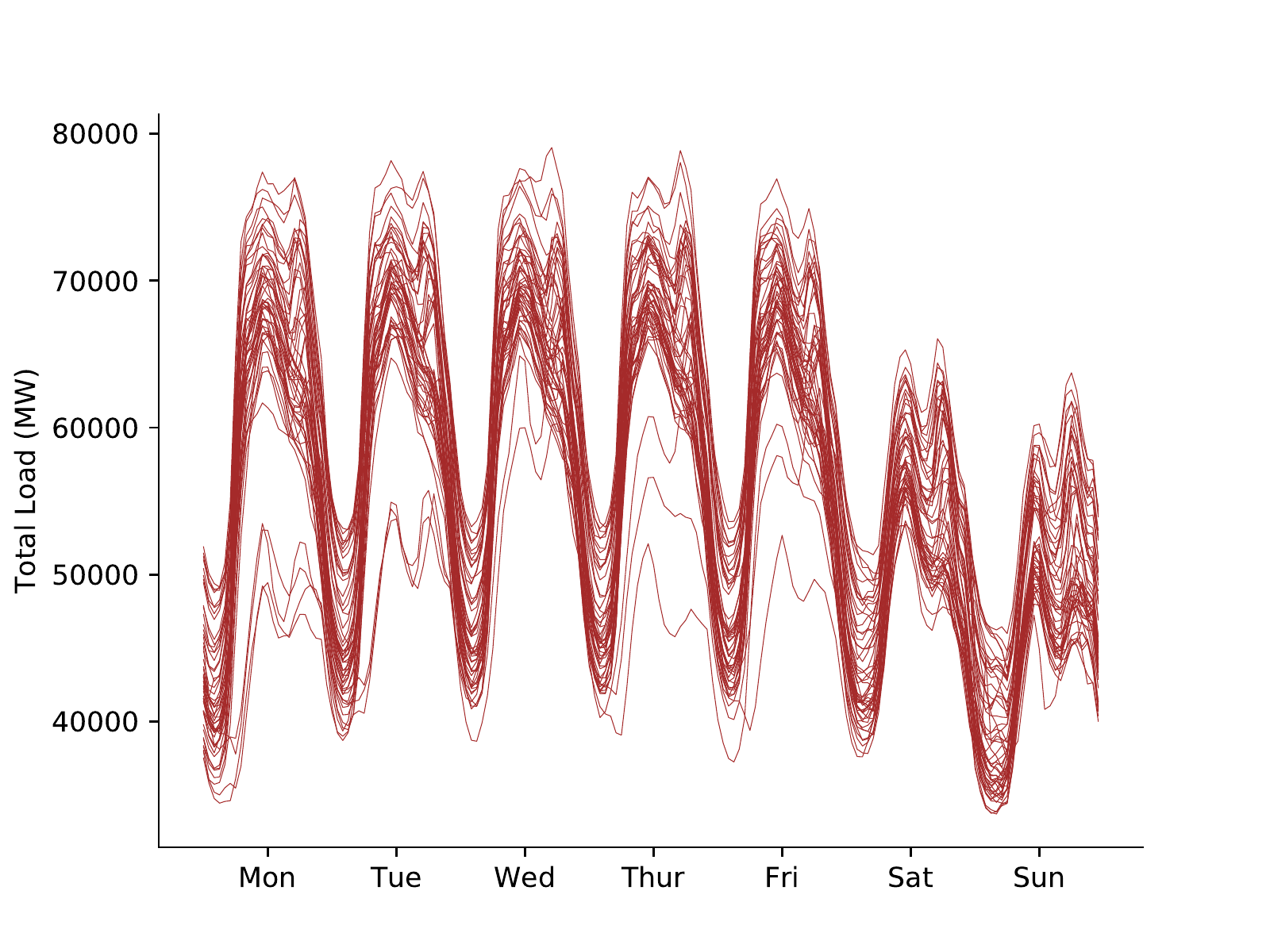}
	\caption{52 trajectories of the data}
	\label{fig:data}
\end{figure}

\paragraph{Brief summary of the data.}
The hourly actual load of electricity data consists of $52$ trajectories representing the $52$ weeks in a year (cf.\ \cref{fig:data}). Each trajectory consists of the 7 days of the week. The data is recorded from $12\colon00$\,a.m.\ to $11\colon59$\,p.m.\ of each day. This means that we are employing a $52\times168$ dimensional data where 168 stages represent the $24\times7$ hours of each week.
The data shows a similar pattern in all the weeks. The actual total load of electricity tends to be higher during the weekdays than on the weekends. It is also higher during the daytime than in the nighttime. The maximum hourly mean of the data is \SI{70290.4}{\MW} and the minimum hourly mean is \SI{38775.3}{\MW}. 
A special characteristic of the data is the presence of outliers which specifically represent holidays in Germany. For example, four of the holidays in the year 2017 fell on Monday (e.g., Easter Monday (April 17\textsuperscript{th}), Labor day (May 1\textsuperscript{st}), Whit day (June 5\textsuperscript{th}), Christmas (Dec 25\textsuperscript{th})). These Monday-outliers are visible in \cref{fig:data}.

\subsection{Approximation by a scenario lattice}
The data in \cref{fig:data} is limited. There are only $52$ trajectories in the data. We use the stochastic approximation procedure to discretize this data (\cref{alg2}). Ideally, to obtain convergence of the stochastic approximation procedure, it is necessary to have an unlimited number of sample paths available. The step size $\alpha_k = \frac{1}{3000+k}$, 
	where $k$ is the stochastic approximation iteration, turned out to be useful for this data.

\begin{figure}[H]
	\centering
	\includegraphics[width=0.75\textwidth]{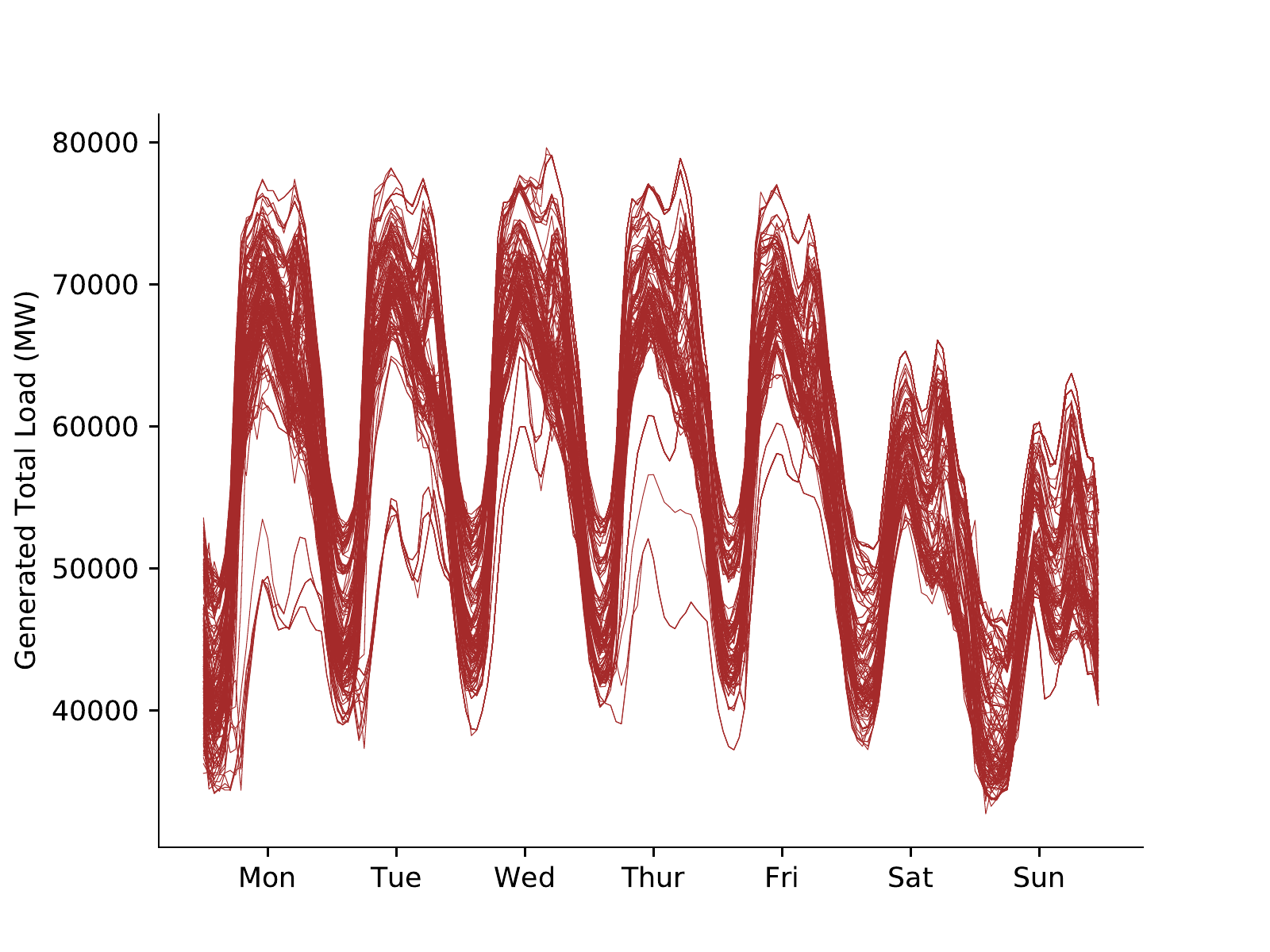}
	\caption{\SI{1000} new trajectories generated from the sample in \cref{fig:data} by employing \cref{alg3}}
	\label{fig:generated}
\end{figure}

Based on this data and using \cref{alg3}, \cref{fig:generated} shows a sample of new and additionally generated scenarios.
These generated trajectories follow the same pattern as in the original data (cf.\ \cref{fig:data}). All the important and necessary characteristics of the original sample is also captured in this data. Therefore, these trajectories totally represent the original sample without any loss of information.

To use 
\href{https://kirui93.github.io/ScenTrees.jl/stable/}{\texttt{ScenTrees.jl}}\textsuperscript{\ref{fn:ScenTrees}}  
we fix the branching structure of the approximating scenario lattice and the number of iterations. In this case, we generate a scenario lattice with 5 nodes at each stage and $2.0\times10^6$ iterations. 
\begin{figure}[H]
	\centering
	\begin{subfigure}[t]{0.49\textwidth}
		\includegraphics[trim=00 18 00 36, clip, width=1.0\textwidth]{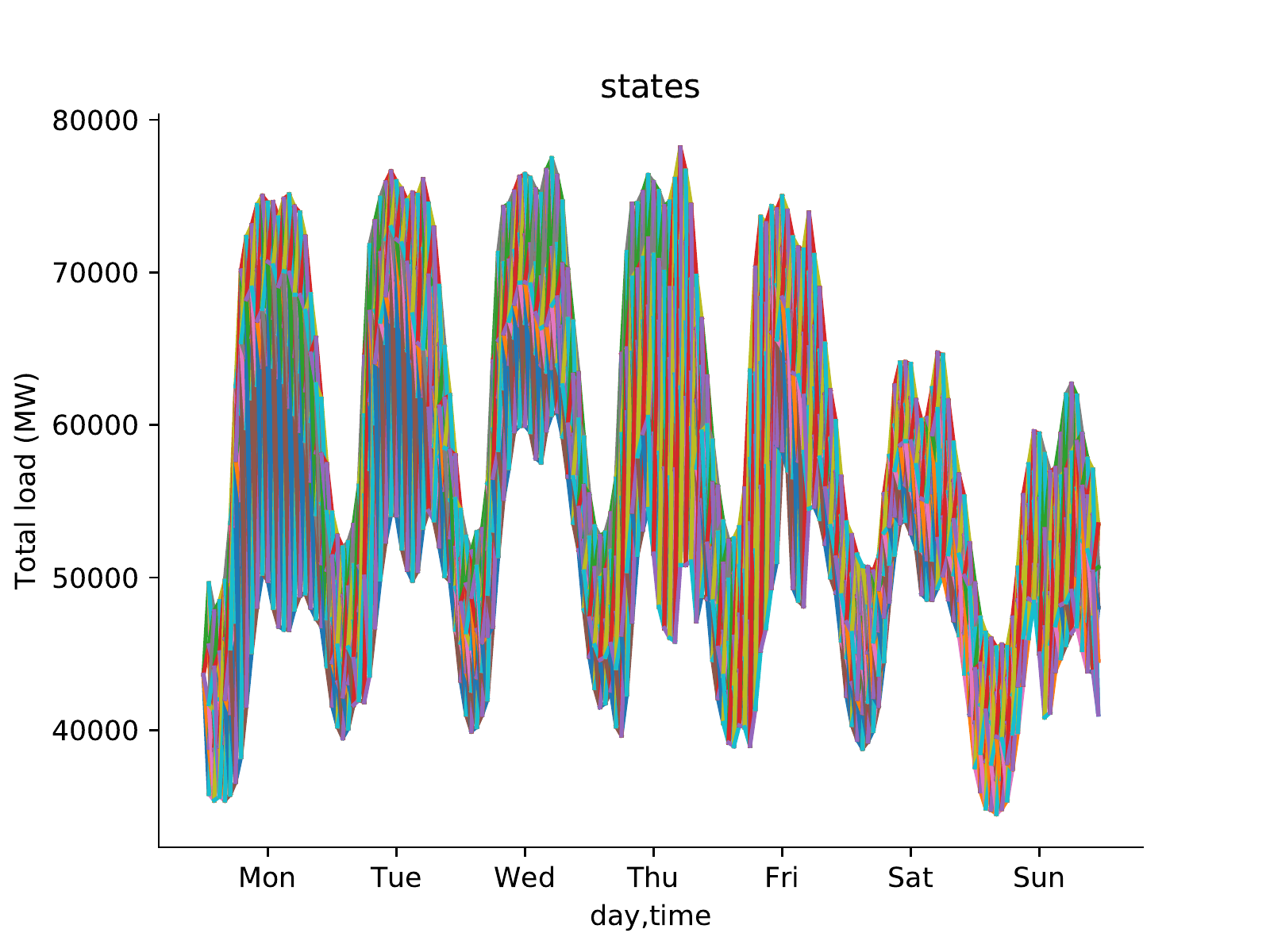}
		\caption{Energy demand during the entire week}
	\end{subfigure}\hfill
	\begin{subfigure}[t]{0.49\textwidth}
		\includegraphics[trim=15 18 55 55, clip, width=1.0\textwidth]{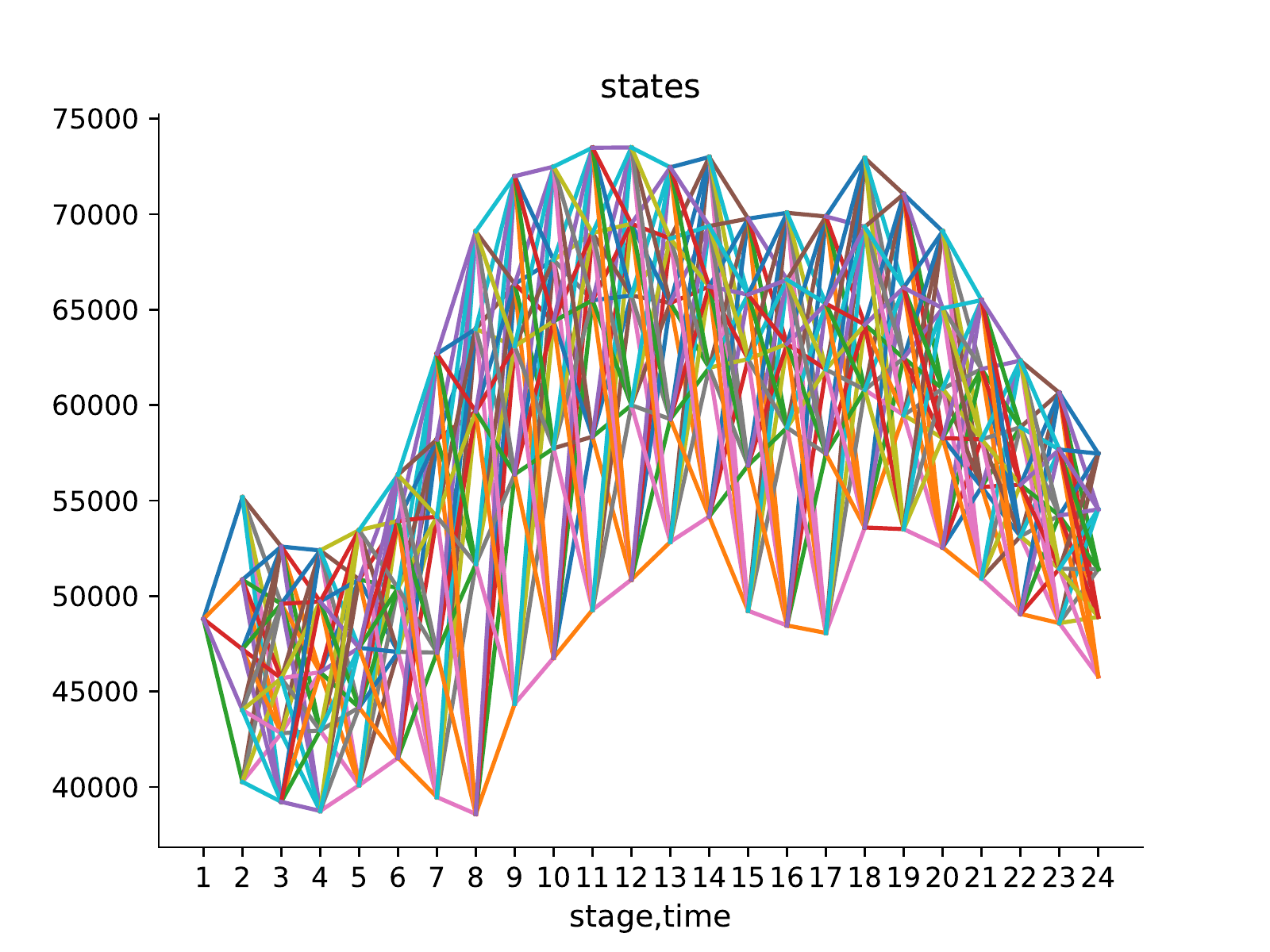}
		\caption{Detailed, hourly view of Saturday}
	\end{subfigure}
	\caption{Scenario lattice generated from \cref{fig:data} using \cref{alg3} and \cref{alg2}}
	\label{fig:resLat}
\end{figure}  
The scenario lattice in \cref{fig:resLat} has $5$ branches connected to each node in each stage. There are a total of $836$ nodes and approximately $5.4\times10^{116}$ scenarios possible in this lattice. 

It is important to note that the generated scenario lattice is able to represent almost all the characteristics in the original data in \cref{fig:data}. This lattice recognizes the patterns in the original data and is able to discretize the data and recover this pattern. To determine the quality of approximation of the scenario lattice, we employ the transportation distance (cf.~\eqref{aad}). Ideally, one would want to find a scenario lattice which minimizes this distance. The generated scenario lattice in \cref{fig:resLat} has a transportation distance of \SI{1,203}{\MW} per stage.

\section{Conclusions}\label{sec:Summary}
In this paper, we present new and fast implementations to generate scenario trees and scenario lattices for decision making under uncertainty.
The package
\begin{center}
	\href{https://kirui93.github.io/ScenTrees.jl/stable/}{\texttt{ScenTrees.jl}}\textsuperscript{\ref{fn:ScenTrees}}
\end{center}
collects efficient and freely available implementations of all algorithms.
We also provided explicit concepts to measure the quality of the approximation by employing distances, particularly the transportation distance. This distance constitutes an essential tool in stochastic optimization. We introduce several techniques to constructively obtain scenario tress from samples, which are observed trajectories. 
We particularly elaborate on methods to find scenario trees and lattices with limited data only. We specifically apply scenario lattice generation methods on a limited electricity load data.

Our results prove that the new algorithms for generating scenario trees and scenario lattices are important as they are able to discretize the data such that all the information in the data is captured. These results also show that our implementation is highly competitive in terms of computational performance.


%

\nocite{Pflug2001}  

\bibliographystyle{abbrvnat}
\bibliography{LiteraturAlois}

\appendix

\section{The Kantorovich/\,Wasserstein distance}\label{sec:Wasserstein}
Recall the definition of the Kantorovich/\,Wasserstein distance $\mathsf{d}_r(P,\tilde{P})$ of order $r$ for two (Borel) random measures $P$ and $\tilde{P}$ on $\mathbb{R}^m$:
\begin{align}
	\mathsf{d}_r(P,\tilde{P}) :=\inf_{\pi} & \left(\iint_{\mathbb{R}^m \otimes \mathbb{R}^m} \| \xi - \tilde{\xi}\|^{r} \, \pi(\mathrm d\xi,\mathrm d\tilde{\xi})\right)^{1/r}\nonumber \\
	\textnormal{s.t. } & \pi\left(A\times \mathbb{R}^m \right)= P(A)\label{eq:12} \\
	& \pi\left(\mathbb{R}^m \times B \right)=\tilde{P}(B) \label{eq:13} .
\end{align}
The infimum is over all Borel measures $\pi$ on $\mathbb{R}^m \times \mathbb{R}^m$ with given marginals $P$ and $\tilde{P}$.
These measures are called \emph{transportation plans}.  If $X$ and $\tilde X$ are $\mathbb R^m$-valued random variables, then their distance is defined as the distance of the corresponding image measures $P^X$ and $P^{\tilde X}$.

The Kantorovich/\,Wasserstein distance does not always have a representation via a transportation map. If, however, $P$ is atomless on a separable space $\Xi$, then it can be employed to approximate any measure in Kantorovich/\,Wasserstein distance arbitrarily close via a transport plan. I.e., for every measure $\tilde P$ on $\tilde\Xi$ and $\varepsilon>0$ there is a transport map $T\colon\Xi\to\tilde\Xi$ such that
\[\mathsf{d}\big(\tilde{P}, P^T\big)<\varepsilon\ \text{ and }\ \mathsf d\big(P^T,P\big)<\mathsf d(\tilde P,P)+\varepsilon.\]

To accept the assertion recall from \citet{Bolley2008} that there are discrete measures $P_n=\sum_{i=1}^n p_i \delta_{x_i}$ and $\tilde P_m=\sum_{j=1}^m \tilde p_j \delta_{y_j}$ such that $\mathsf d(P,P_n)<\varepsilon/2$ and $\mathsf d(\tilde P,\tilde P_m)<\varepsilon/2$. Further, there is a Dirichlet tessellation of 
$\Xi$ in disjoint sets, $\Xi=\dot \bigcup_{i=1}^n \Xi_i$, such that
\[\Xi_i\subset\{\xi\colon d(\xi,x_i)\le d(\xi,x_{i^\prime})\text{ for all }i^\prime=1,\dots,n \}\] and without loss of generality we may assume that $p_i=P(\Xi_i)$. 

Let $\pi$ denote the transport plan so that $\mathsf d(P_n,\tilde P_m)=\sum_{i=1}^n\sum_{j=1}^m d(x_i,y_j)\pi_{i,j}$, satisfying the marginal constraints~\eqref{eq:12} and~\eqref{eq:13}, i.e., $p_i=\sum_{j=1}^m \pi_{i,j}$ and $\tilde p_j=\sum_{i=1}^n \pi_{i,j}$. As $P$ does not give masses to atoms, one may dissect the sets $\Xi_i=\dot\bigcup_{j=1}^m \Xi_{i,j}$ so that $P(\Xi_{i,j})=\pi_{i,j}$. It holds that $\Xi=\dot\bigcup_{i=1,j=1}^{n,\ \ m}\Xi_{i,j}$, so that the transport map
\[T(x):= y_j\text{ if }x\in\Xi_{i,j}\]
is well defined. By construction, $P^T=\tilde P_m$, so that $\mathsf d(\tilde P,P^T)<\varepsilon$ and $\mathsf d(P,P^T)<\mathsf d(P,\tilde P)+\varepsilon$. 
\[\mathsf d(\tilde P,P^T)\le \mathsf d(\tilde P,\tilde P_m),\]
the assertion.

\section{The nested distance}\label{sec:NestedA}
The nested distance introduced in \cite{Pflug2009, PflugPichler2011} generalizes the Kantoro\-vich/\,Was\-ser\-stein distance. Its definition given in~\eqref{eq:1}--\eqref{eq:4} above is for the laws of stochastic processes, the image measures on filtered probability spaces. 

The following mathematical definition of the nested distance generalizes the recursive description above. It is an obvious generalization of the Wasserstein distance, which  \Cref{sec:Wasserstein} elaborates in detail.
The nested distance, denoted $\dd$, is
\begin{align}
	\inf_\pi\E_\pi d(\xi,\eta)&=\iint_{\Xi^T\times\Xi^T}d(\xi,\eta)\,\pi(\mathrm d \xi,\mathrm d\eta)\label{eq:1}\\
	\text{subject to }& \pi(A\times \Xi^T\mid \mathcal F_t\times \mathcal F_t)=P(A\mid \mathcal F_t),\quad A\in\mathcal F_T \text{ and }\label{eq:3}\\
	& \pi(\Xi^T\times B\mid\mathcal F_t\times\mathcal F_t)=Q(B\mid\mathcal F_t), \quad B\in\mathcal F_T\label{eq:4}\\
	&\qquad\qquad\text{ for all stages }t\in\{1,\dots,T\},\nonumber
\end{align}
where the infimum is among all bivariate probability measures $\pi$ satisfying the marginal constraints~\eqref{eq:3} and~\eqref{eq:4} at all stages $t\in\{1,\dots,T\}$.

We note that for random variables, or stochastic processes in one stage, the definition of the nested distance~\eqref{eq:1}--\eqref{eq:4} coincides with the Wasserstein or Kantorovich distance. Both are transportation distances, but the nested distance additionally takes the evolving information into account via~\eqref{eq:3} and~\eqref{eq:4}.

The construction principle of the nested distance is (backwards) recursive. To compute the nested distance of two trees of same height we isolate the roots from the tree and consider the root, and all remaining subtrees separately:
\begin{enumerate}[label=(\roman*)]
	\item \label{enu:1} The root nodes of both trees have a specified distance, which is available by employing the distance $d_t$ at the stage $t$.
	\item \label{enu:2} Each successor node of the root node is the root node of a subtree. Further, a transition probability from the root node to each subtree is available by employing the transportation distance, applied do the subtrees.

	Each combination of subtrees from both trees has a nested distance, which is available by recursion. These nested distances of subtrees constitute (abstract) costs. The transportation distance, as specified above, makes a distance of the successor nodes (i.e., subtrees) available.
\end{enumerate}

The nodes at the final stage $t=T$ do not have subtrees, such that the component~\ref{enu:2} does not apply at stage $t=T$: by~\ref{enu:1}, the nested distance is $d_T$ at stage $t=T$.
The nested distance of intermediate trees at the same stage is given by recombining the roots and the subtrees, i.e., the nested distance is \ref{enu:1}\,+\,\ref{enu:2}.

The nested distance finally is the distance of the trees from their root node $0$ and available in a backwards recursive way.

To interpret this definition, the nested distance between two multistage probability distributions is obtained by minimizing over all transportation plans $\pi$ transporting one distribution into the other, which are compatible with the filtration structures in addition. For a single period, i.e. $T=1$, the nested distance coincides with the Kantorovich/\,Wasser\-stein distance. The nested distance is crucial for multistage stochastic optimization by the following theorem on stability of stochastic multistage optimization problems (cf.\ \cite{PflugPichler2011}).

\begin{theorem}\label{thm:SD}
Let $\mathbb{P}:=(\Xi, (\mathcal F_t)_{t=1}^T, P)$ ($\tilde{\mathbb{P}}$, resp.) be a filtered probability space. Consider the multistage stochastic optimization problem
\[
	v(\mathbb{P}):=\inf\left\{ \E_\xi^{\mathbb{P}} \big(Q(\xi,x)\big)\colon x\lhd\Fc \right\},
\]
where $Q$ is convex in the decisions $x=(x_1,\dots, x_T)$ for any $\xi$ fixed, and Lipschitz with
constant $L$ in the scenario process $\xi=(\xi_1, \dots, \xi_T)$ for any $x$ fixed. The constraint $x\lhd \Fc$ means that the decisions can be random variables, but must be adapted to the filtration $\Fc$, i.e., must be nonanticipative. Then the objective values $v(\mathbb{P})$ and $v(\tilde{\mathbb{P}})$ calculated with the two different probability models for the scenario process $\mathbb{P}$ and $\tilde{\mathbb{P}}$, satisfy
\[
	\left|v(\mathbb{P})-v(\tilde{\mathbb{P}})\right|\le L\cdot\nestr(\mathbb{P},\tilde{\mathbb{P}}) \, .
\]
\end{theorem}

\vfill\hfill\scriptsize\today,\ \jobname.pdf
\end{document}